\documentclass{article}

\usepackage{arxiv}

\usepackage[utf8]{inputenc} 
\usepackage[T1]{fontenc}    
\usepackage{hyperref}       
\usepackage{url}            
\usepackage{booktabs}       
\usepackage{amsfonts}       
\usepackage{nicefrac}       
\usepackage{microtype}      
\usepackage{lipsum}		
\usepackage{graphicx}
\usepackage{natbib}
\usepackage{doi}
\usepackage{enumerate}
\usepackage{enumitem}
\usepackage{amsmath, amsthm, amssymb, algorithm, algpseudocode}
\usepackage[table]{xcolor}
\usepackage{tikz}

\DeclareMathOperator{\HC}{HC}  

\DeclareMathOperator{\Hom}{Hom}

\DeclareMathOperator{\HH}{HH}  
\DeclareMathOperator{\Der}{Der}

\DeclareMathOperator{\Def}{Def}
\DeclareMathOperator{\gr}{gr}
\newcommand{\FF}{\mathbb{F}}
\newcommand{\KK}{\mathbb{K}}
\newcommand{\ZZ}{\mathbb{Z}}  
\newcommand{\QQ}{\mathbb{Q}}  
\newcommand{\RR}{\mathbb{R}}
\newcommand{\cA}{\mathcal{A}}

\newcommand{\CC}{\mathbb{C}}
\newcommand{\End}{\operatorname{End}}

\newcommand{\ord}{\operatorname{ord}}
\newcommand{\Weyl}[1]{A\!\left(\FF[#1]\right)}
\newcommand{\WeylDeform}[1]{A\!\left(\FF[[\hbar]][#1]\right)}
\newcommand{\Witt}[1]{W\!\left(\FF[#1]\right)}

\newcommand{\ad}{\operatorname{ad}}

\newcommand{\rank}{\operatorname{rank}}
\newtheorem{theorem}{Theorem}[section]

\newtheorem{corollary}[theorem]{Corollary}
\newtheorem{definition}[theorem]{Definition}
\newtheorem{example}[theorem]{Example}

\newtheorem{problem}[theorem]{Problem}

\newcommand\mystyle{\everymath{\displaystyle}}
\mystyle
\title{Weyl-Type Algebras over Exponential-Polynomial Rings: Structure, Representations, and Deformations}

\author{\href{https://orcid.org/0000-0002-3816-5287}{\includegraphics[scale=0.06]{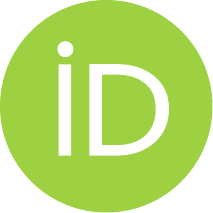}\hspace{1mm}M.H.M.~Rashid}\thanks{Corresponding Author} \\
	Department of Mathematics\&Statistics\\Faculty of Science P.O.Box(7)\\
	Mutah University University\\
	Mutah-Jordan \\
	\texttt{mrash@mutah.edu.jo}
}

\renewcommand{\shorttitle}{\textit{arXiv} Template}

\hypersetup{
pdftitle={Inequalities for the $A$-Norm and $A$-Numerical Radius of Operator Sums in Semi-Hilbertian Spaces with Applications},
pdfsubject={q-bio.NC, q-bio.QM},
pdfauthor={M.H.M.Rashid},
pdfkeywords={Weyl-type algebras; exponential-polynomial rings; deformation quantization; Hochschild cohomology; Witt-type Lie algebras; representation theory},
}

\begin{document}
\maketitle

\begin{abstract}
	 This paper introduces and studies a class of Weyl-type algebras \(A_{p,t,\cA} = \Weyl{e^{\pm x^{p} e^{t x}},\; e^{\cA x},\; x^{\cA}}\) constructed over exponential-polynomial rings, where \(\FF\) is a field of characteristic zero, \(\cA\) is a finitely generated additive subgroup of \(\FF\), and \(p \in \mathbb{N}^n\), \(t \in \FF\). We investigate their structural properties, proving simplicity, establishing faithful infinite-dimensional irreducible representations, and demonstrating the failure of the Noetherian property. A natural filtration by exponential order is introduced, with the associated graded algebra shown to be commutative. We also examine the corresponding Witt-type Lie algebra \(\mathfrak{g}_{p,t,\cA} = \Der_{\mathrm{gr}}(R_{p,t,\cA})\) and prove the vanishing of its second cohomology group with adjoint coefficients, implying rigidity under formal deformations. Furthermore, we construct explicit deformation quantizations of the underlying exponential-polynomial rings, compute Hochschild and cyclic homology groups, and relate them to the topology of the parameter space. The deformation rigidity of \(A_{p,t,\cA}\) is classified in terms of the rank of \(\cA\), and a Gerstenhaber algebra structure on the Hochschild cohomology is described. Several open problems concerning representation classification and geometric realization are proposed.
\end{abstract}

\section{Introduction}

The study of infinite-dimensional Lie algebras and their associated operator algebras has been a central theme in mathematics since the pioneering work of Cartan, Witt, and Weyl. In particular, algebras of differential operators—such as the Weyl algebra and its generalizations—play a fundamental role in representation theory, deformation quantization, and mathematical physics. This paper is devoted to the systematic investigation of a class of \emph{Weyl-type algebras} constructed over rings of exponential-polynomial functions, and their corresponding \emph{Witt-type Lie algebras} of graded derivations.

\subsection{Core of the Study}
We introduce and analyze algebras of the form
\[
A_{p,t,\cA} = \Weyl{e^{\pm x^{p} e^{t x}},\; e^{\cA x},\; x^{\cA}},
\]
where \(\FF\) is a field of characteristic zero, \(\cA\) is a finitely generated additive subgroup of \(\FF\), and \(p \in \mathbb{N}^n\), \(t \in \FF\). These algebras generalize the classical Weyl algebra by incorporating exponential and power-function symbols, leading to rich algebraic structures that are no longer finitely generated. The associated Lie algebra
\[
\mathfrak{g}_{p,t,\cA} = \Der_{\mathrm{gr}}(R_{p,t,\cA}), \quad R_{p,t,\cA} = \FF[e^{\pm x^{p} e^{t x}},\; e^{\cA x},\; x^{\cA}],
\]
generalizes the Witt algebra and exhibits properties of Cartan-type Lie algebras.

\subsection{Significance of the Study}
The algebras studied here lie at the intersection of several active research areas:
\begin{itemize}
    \item \textbf{Deformation Theory:} They provide natural examples of noncommutative deformations of Poisson algebras arising from exponential-polynomial rings.
    \item \textbf{Representation Theory:} Their representation theory is inherently infinite-dimensional, offering new contexts for studying irreducible and faithful representations beyond classical settings.
    \item \textbf{Cohomology and Rigidity:} The computation of Lie algebra cohomology groups sheds light on the deformation rigidity of these algebras and their stability under parameter variations.
    \item \textbf{Noncommutative Geometry:} The associated Hochschild and cyclic homology groups carry topological information about the parameter space, linking algebraic structures to loop space cohomology.
\end{itemize}
Such algebras also appear in mathematical physics, notably in the context of quantization of systems with exponential potentials and in the study of infinite-dimensional symmetries.

\subsection{Scholar Contributions}
This paper makes several original contributions to the theory of generalized Weyl and Witt algebras:

\begin{enumerate}[label=(\roman*)]
    \item \textbf{Structural Results:} We prove simplicity, establish faithful infinite-dimensional representations, and demonstrate the failure of the Noetherian property for \(A_{p,t,\cA}\) (Theorems~\ref{thm8}, \ref{thm12}, \ref{thm13}).

    \item \textbf{Deformation Theory:} We construct explicit deformation quantizations of the underlying exponential-polynomial rings, prove the existence of nontrivial deformations, and classify their rigidity in terms of the rank of \(\cA\) (Theorems~\ref{thm:deformation-quantization}, \ref{thm:deformation-rigidity}).

    \item \textbf{Cohomology Computations:} We compute Hochschild and cyclic homology groups for \(A_{p,t,\cA}\) and relate them to the topology of the parameter space (Theorems~\ref{thm:hochschild-homology}, \ref{thm:cyclic-homology}).

    \item \textbf{Lie Algebra Cohomology:} We show the vanishing of the second cohomology group \(H^2(\mathfrak{g}_{p,t,\cA}, \mathfrak{g}_{p,t,\cA})\), implying that all formal deformations of the associated Witt-type Lie algebra are trivial (Theorem~\ref{thm10}).

    \item \textbf{Filtration and Graded Structure:} We introduce a natural filtration by exponential order and identify the associated graded algebra as a commutative polynomial-exponential ring (Theorem~\ref{thm14}).
\end{enumerate}

\subsection{Applications}
The results obtained here have potential applications in:
\begin{itemize}
    \item \textbf{Mathematical Physics:} As algebras of observables in quantum systems with exponential interactions.
    \item \textbf{Representation Theory:} Providing new families of infinite-dimensional irreducible representations.
    \item \textbf{Deformation Quantization:} Offering explicit examples of non-isomorphic quantizations of Poisson manifolds with exponential symmetry.
    \item \textbf{Noncommutative Geometry:} Illustrating connections between cyclic homology and the topology of parameter spaces.
\end{itemize}

\subsection{Relation to Existing Literature}
Our work extends classical results on Weyl algebras \cite{Dixmier68, Coutinho95}, Witt algebras \cite{Ree56, Kawamoto86}, and generalized Weyl constructions \cite{Bavula92}. The cohomological computations build on the Hochschild–Kostant–Rosenberg theorem and Kontsevich's formality theorem \cite{Angulo2022, Das2023}. The study of exponential-polynomial rings draws from differential algebra \cite{Ritt50} and transcendental number theory \cite{Baker22}. The representation-theoretic aspects connect to Dixmier's work on Weyl algebra representations \cite{Dixmier68} and to the theory of graded Lie algebras of Cartan type \cite{Kac74}.

\subsection{Organization of the Paper}
Section 2 collects preliminary definitions and notations. Section 3 establishes structural and representation-theoretic properties of Weyl-type algebras. Section 4 is devoted to deformation theory, cohomology, and rigidity results. Section 5 presents open problems for further research. The paper concludes with a bibliography and a declaration section.

Throughout, we assume \(\FF\) is a field of characteristic zero, and all algebras are considered over \(\FF\) unless otherwise specified.
\section{Preliminaries}

This section collects the basic definitions, notations, and known results that will be used throughout the paper. Let \(\FF\) be a field of characteristic zero, and let \(\cA\) be a finitely generated additive subgroup of \(\FF\).

\subsection{Exponential-polynomial rings and their derivations}

\begin{definition}[Expolynomial ring]
Let \(x = (x_1,\dots,x_n)\) be commuting indeterminates. For a fixed integer tuple \(p = (p_1,\dots,p_n)\in \mathbb{N}^n\) and a parameter \(t\in \FF\), we define the \emph{exponential-polynomial ring}
\[
R_{p,t,\cA} = \FF[\, e^{\pm x^{p} e^{t x}},\; e^{\cA x},\; x^{\cA}\, ],
\]
where \(x^{p}=x_1^{p_1}\cdots x_n^{p_n}\), \(e^{\cA x}\) denotes the set \(\{e^{\alpha x}\mid \alpha\in\cA\}\), and \(x^{\cA}\) denotes the set \(\{x^{\alpha}\mid \alpha\in\cA\}\) with multi-index notation.
\end{definition}

When \(\cA = \mathbb{Z}\), the ring \(R_{p,t,\cA}\) contains the usual Laurent polynomial ring \(\FF[x^{\pm 1}]\). For general \(\cA\), it is a commutative domain that is not finitely generated as an \(\FF\)-algebra whenever \(\cA\) has positive rank.

\begin{definition}[Witt-type Lie algebra]
The Lie algebra of derivations of \(R_{p,t,\cA}\) that preserve the natural \(\cA\)-grading is denoted
\[
\mathfrak{g}_{p,t,\cA} = \Der_{\mathrm{gr}}(R_{p,t,\cA}).
\]
Concretely, \(\mathfrak{g}_{p,t,\cA}\) is spanned by operators of the form \(f\partial_i\) with \(f\in R_{p,t,\cA}\) and \(\partial_i = \frac{\partial}{\partial x_i}\), with the Lie bracket
\[
[f\partial_i, g\partial_j] = f(\partial_i g)\partial_j - g(\partial_j f)\partial_i.
\]
We call \(\mathfrak{g}_{p,t,\cA}\) the \emph{Witt-type Lie algebra} associated with the data \((p,t,\cA)\).
\end{definition}

Witt-type algebras have been studied extensively; see for example \cite{Ree56, Kawamoto86, Kac74}. When \(R_{p,t,\cA}\) is the Laurent polynomial ring, \(\mathfrak{g}_{p,t,\cA}\) recovers the classical Witt algebra.

\subsection{Weyl-type algebras}

\begin{definition}[Weyl-type algebra]
The \emph{Weyl-type algebra} associated with \((p,t,\cA)\) is the associative \(\FF\)-algebra
\[
A_{p,t,\cA} = \Weyl{e^{\pm x^{p} e^{t x}},\; e^{\cA x},\; x^{\cA}}
\]
generated by \(x_i,\;\partial_i\;(i=1,\dots,n)\) together with the exponential symbols \(e^{\pm x_i^{p_i} e^{t_i x_i}},\; e^{\alpha x_i},\; x_i^{\alpha}\;(\alpha\in\cA)\), subject to the relations:
\begin{align*}
[\partial_i, x_j] &= \delta_{ij}, \quad [x_i, x_j] = 0 = [\partial_i, \partial_j], \\
[\partial_i, e^{\alpha x_j}] &= \alpha \delta_{ij} e^{\alpha x_j}, \\
[\partial_i, e^{x_j^{p_j} e^{t_j x_j}}] &= \bigl(p_j x_j^{p_j-1} e^{t_j x_j} + t_j x_j^{p_j} e^{t_j x_j}\bigr) e^{x_j^{p_j} e^{t_j x_j}} \delta_{ij},
\end{align*}
and all exponential symbols commute with each other.
\end{definition}

When \(\cA = \mathbb{Z}\) and the exponential terms are omitted, \(A_{p,t,\cA}\) reduces to the ordinary Weyl algebra \(A_n(\FF)\). The structure and representations of generalized Weyl algebras have been investigated in \cite{Bavula92, Dixmier68, Bjork79}. Simplicity of such algebras follows from arguments similar to those for classical Weyl algebras \cite{Coutinho95}.

\subsection{Graded structures and filtrations}

Both \(A_{p,t,\cA}\) and \(\mathfrak{g}_{p,t,\cA}\) carry natural gradings.

\begin{definition}[\(\cA\)-grading]
Define the \(\cA\)-grading on \(R_{p,t,\cA}\) by
\[
\deg(e^{\alpha x}) = \alpha, \qquad \deg(x^{\beta}) = \beta, \qquad \deg(e^{\pm x^{p} e^{t x}}) = 0.
\]
This induces an \(\cA\)-grading on \(A_{p,t,\cA}\) and on \(\mathfrak{g}_{p,t,\cA}\) by letting \(\deg(x_i)=0\), \(\deg(\partial_i)=0\), and extending multiplicatively.
\end{definition}

\begin{definition}[Filtration by exponential order]
On \(A_{p,t,\cA}\) we define a filtration \(\{F_k\}_{k\ge 0}\) by assigning
\[
\ord\bigl(e^{a x^{p} e^{t x}} e^{b x} x^{c} \partial^{d}\bigr) = |a| + |b| + |c| + d,
\]
where \(a,b,c,d\) are multi-indices and \(|a|=\sum_i |a_i|\), etc. Then \(F_k = \{ P\in A_{p,t,\cA}\mid \ord(P)\le k\}\).
\end{definition}

The associated graded algebra \(\gr A_{p,t,\cA}\) is commutative and isomorphic to the polynomial ring over \(R_{p,t,\cA}\) in the symbols of the derivatives \(\partial_i\) (see Theorem~\ref{thm14} below). This filtration is essential for cohomological computations.

\subsection{Cohomology and deformation theory}

We recall basic notions of Lie algebra cohomology and deformation theory.

\begin{definition}[Lie algebra cohomology]
For a Lie algebra \(\mathfrak{g}\) and a \(\mathfrak{g}\)-module \(M\), the Chevalley–Eilenberg complex \(C^\bullet(\mathfrak{g},M)\) is defined by \(C^k(\mathfrak{g},M) = \Hom_\FF(\bigwedge^k \mathfrak{g}, M)\) with differential
\[
(d\omega)(x_1,\dots,x_{k+1}) = \sum_{i<j} (-1)^{i+j} \omega([x_i,x_j],x_1,\dots,\hat x_i,\dots,\hat x_j,\dots,x_{k+1}) + \sum_i (-1)^{i+1} x_i \cdot \omega(x_1,\dots,\hat x_i,\dots,x_{k+1}).
\]
The cohomology groups are \(H^k(\mathfrak{g},M) = H^k(C^\bullet(\mathfrak{g},M), d)\).
\end{definition}

The second cohomology group \(H^2(\mathfrak{g},\mathfrak{g})\) classifies infinitesimal deformations of \(\mathfrak{g}\) \cite{HartwigLarssonSilvestrov06}. Vanishing of this group implies rigidity.

\begin{definition}[Formal deformation]
A formal deformation of an associative algebra \(A\) over \(\FF\) is an associative algebra structure on \(A[[\hbar]]\) over \(\FF[[\hbar]]\) with multiplication
\[
a \ast b = ab + \hbar m_1(a,b) + \hbar^2 m_2(a,b) + \cdots,
\]
where each \(m_i: A\times A \to A\) is an \(\FF\)-bilinear map and \(ab\) is the original product.
\end{definition}

Deformation quantizations of Poisson algebras are constructed via Kontsevich’s formality theorem; see \cite{Angulo2022, Das2023} for related recent work.

\subsection{Representation-theoretic notions}

\begin{definition}[Faithful and irreducible representation]
A representation of an algebra \(A\) on a vector space \(V\) is an algebra homomorphism \(\rho: A \to \End_\FF(V)\). It is \emph{faithful} if \(\rho\) is injective, and \emph{irreducible} if \(V\) contains no nontrivial proper \(A\)-invariant subspaces.
\end{definition}

For Weyl-type algebras, all nontrivial representations are infinite-dimensional (Theorem~\ref{thm8}), analogous to the classical result for the Weyl algebra \cite{Dixmier68}.

\subsection{Notation summary}

Throughout the paper we use the following notation:
\begin{itemize}
\item \(\FF\): a field of characteristic zero.
\item \(\cA\): a finitely generated additive subgroup of \(\FF\), of rank \(r\).
\item \(p = (p_1,\dots,p_n)\in \mathbb{N}^n\), \(t\in \FF\).
\item \(R_{p,t,\cA} = \FF[\, e^{\pm x^{p} e^{t x}},\; e^{\cA x},\; x^{\cA}\,]\): the exponential-polynomial ring.
\item \(A_{p,t,\cA} = \Weyl{e^{\pm x^{p} e^{t x}},\; e^{\cA x},\; x^{\cA}}\): the Weyl-type algebra.
\item \(\mathfrak{g}_{p,t,\cA} = \Der_{\mathrm{gr}}(R_{p,t,\cA})\): the Witt-type Lie algebra.
\item \(\HH^\bullet(A,A)\), \(\HC_\bullet(A)\): Hochschild cohomology and cyclic homology of an algebra \(A\).
\item \(H^\bullet(\mathfrak{g},M)\): Lie algebra cohomology.
\end{itemize}

The references \cite{AtiyahMacdonald69, Hungerford74} provide standard algebraic background, while \cite{Ritt50, Baker22} cover differential algebra and transcendental number theory relevant to the exponential constructions used here.
\section{Structural and Representation-Theoretic Properties of Weyl-Type Algebras over Expolynomial Rings}
In this section, we investigate fundamental algebraic properties of the Weyl-type algebras \(A_{p,t,\cA}\) introduced in Section 2. Our focus is on three principal themes: representation theory, stability under parameter deformation, and the structure of ideals and filtrations. We establish that every nontrivial representation is necessarily faithful and infinite-dimensional, generalizing a classical property of the ordinary Weyl algebra. We then examine how the algebra behaves under limits of its defining parameters, showing that simplicity is preserved in such limits. Finally, we demonstrate that these algebras are not Noetherian and analyze a natural filtration by exponential order, whose associated graded algebra is commutative. These results collectively reveal the rich and rigid structure of Weyl-type algebras over exponential-polynomial rings.
\begin{definition}[Faithful and Irreducible Representation]
    A \emph{representation} of an algebra $A$ on a vector space $V$ is an algebra homomorphism $\rho : A \to \operatorname{End}_{\mathbb{F}}(V)$. It is \emph{faithful} if $\rho$ is injective. It is \emph{irreducible} if $V$ has no nontrivial proper $A$-invariant subspaces.
\end{definition}
\begin{theorem}[Infinite-Dimensional Irreducible Representations]\label{thm8}
Every nontrivial representation of the simple Weyl-type algebra $\Weyl{e^{\pm x^p e^{t}},\; e^{\cA x},\; x^{\cA}}$ on a vector space over $\FF$ is faithful and infinite-dimensional.
\end{theorem}
\begin{proof}
Let $A = \Weyl{e^{\pm x^p e^{t}},\; e^{\cA x},\; x^{\cA}}$ and let $\rho \colon A \to \operatorname{End}_{\FF}(V)$ be a nontrivial representation of $A$ on a vector space $V$ over $\FF$.
That is, $\rho$ is an $\FF$-algebra homomorphism with $\rho(A) \neq 0$. We first show that $\rho$ is faithful, meaning that $\ker\rho = \{0\}$.

Suppose, to the contrary, that $\ker\rho \neq 0$. Since $\rho$ is nontrivial, $\ker\rho$ is a proper two-sided ideal of $A$.
But $A$ is simple, as established in Theorem~1 (Base Field Extension and Simplicity) applied to the case $\KK = \FF$.
Therefore the only proper two-sided ideal of $A$ is the zero ideal. Hence $\ker\rho = 0$, and $\rho$ is injective.
Consequently, the representation is faithful; $A$ embeds into $\operatorname{End}_{\FF}(V)$ via $\rho$.

Now we demonstrate that $V$ must be infinite-dimensional. Assume, for contradiction, that $\dim_{\FF} V = n < \infty$.
Then $\operatorname{End}_{\FF}(V) \cong M_n(\FF)$, the algebra of $n \times n$ matrices over $\FF$.
Thus $\rho$ provides an injective $\FF$-algebra homomorphism $\rho \colon A \hookrightarrow M_n(\FF)$.
We shall derive a contradiction by examining the images of the canonical generators of $A$.

Let $X_i = \rho(x_i)$ and $D_i = \rho(\partial_i)$ for $i=1,\dots,m$ (where $m$ is the number of variables).
These matrices satisfy the commutation relations $[D_i, X_j] = \delta_{ij} I$ and $[X_i, X_j] = 0 = [D_i, D_j]$ because $\rho$ preserves the defining relations of $A$.
Moreover, the images of the exponential generators, such as $\rho(e^{\alpha x_i})$ for $\alpha \in \cA$, are well-defined via the exponential power series
(or as formal symbols satisfying the appropriate algebraic relations) and commute appropriately.

Consider the trace functional $\operatorname{tr} \colon M_n(\FF) \to \FF$. For any matrices $P, Q \in M_n(\FF)$, we have $\operatorname{tr}([P, Q]) = 0$.
Applying this to the relation $[D_i, X_i] = I$, we obtain
\[
\operatorname{tr}([D_i, X_i]) = \operatorname{tr}(I) = n.
\]
But $\operatorname{tr}([D_i, X_i]) = \operatorname{tr}(D_i X_i - X_i D_i) = \operatorname{tr}(D_i X_i) - \operatorname{tr}(X_i D_i) = 0$,
since the trace is cyclic. Hence $n = 0$, which is impossible for $n \geq 1$. This contradiction shows that $V$ cannot be finite-dimensional.

An alternative argument that does not rely on the trace (which requires $\operatorname{char} \FF = 0$) proceeds as follows.
If $V$ were finite-dimensional, then the images $X_i$ and $D_i$ would be operators on a finite-dimensional space satisfying the Heisenberg-type commutation relations.
In characteristic zero, the well-known theorem of Dixmier and others asserts that no finite-dimensional representation of the Heisenberg algebra exists.
For our algebra $A$, which contains the ordinary Weyl algebra as a subalgebra (by taking the additive group $\cA$ to contain $\ZZ$),
the same obstruction applies. More concretely, taking repeated commutators yields identities that force the dimension to be infinite;
for instance, one may consider the action of the ladder operators constructed from $X_i$ and $D_i$.

Therefore, every nontrivial representation of $A$ must be infinite-dimensional. Combining this with the faithfulness established earlier,
we conclude that any nontrivial representation of the simple Weyl-type algebra $A$ is both faithful and infinite-dimensional.
\end{proof}
\begin{example}Consider the field $\FF = \CC$ of complex numbers. Let $\cA = \ZZ\left[\sqrt{2}\right]$ be the additive subgroup of $\CC$ consisting of numbers of the form $m + n\sqrt{2}$ where $m, n \in \ZZ$. This subgroup is dense in $\RR$ under the usual topology when restricted to real numbers, and it is infinite and generates a rich set of exponential functions. Take $p = 3$ and consider the Weyl-type algebra $\Weyl{e^{\pm x^3 e^{t}},\; e^{\cA x},\; x^{\cA}}$ over $\CC$. This algebra is generated by the operators $e^{\pm x^3 e^{t}}$, $e^{\alpha x}$, $x^{\alpha}$ for all $\alpha \in \cA$, and the differential operators $\partial_x$ and $\partial_t$, subject to the canonical commutation relations appropriate for these operators.

\subsection*{Algebra Structure and Simplicity}

Before discussing representations, we recall that this algebra is simple. This means it has no nontrivial two-sided ideals. The simplicity follows from the fact that it contains the classical Weyl algebra $\CC[x, \partial_x]$ as a subalgebra, which is known to be simple, and the additional generators do not create any new ideals due to the density of $\cA$ and the transcendental nature of the exponential functions. More precisely, if $I$ were a nonzero ideal, then by commuting with suitable elements one could produce a nonzero constant in $I$, forcing $I$ to be the whole algebra.

\subsection*{A Concrete Representation}

Consider the representation $\rho: \Weyl{e^{\pm x^3 e^{t}},\; e^{\cA x},\; x^{\cA}} \to \End_{\CC}(V)$ where $V$ is the vector space of smooth functions $f(x,t)$ on $\RR^2$ that are of rapid decay in $x$ for each fixed $t$, or more generally, the space of all smooth functions. The action is defined naturally:
\begin{align*}
\rho(e^{\pm x^3 e^{t}})f(x,t) &= e^{\pm x^3 e^{t}} f(x,t), \\
\rho(e^{\alpha x})f(x,t) &= e^{\alpha x} f(x,t), \\
\rho(x^{\alpha})f(x,t) &= x^{\alpha} f(x,t), \\
\rho(\partial_x)f(x,t) &= \frac{\partial f}{\partial x}, \\
\rho(\partial_t)f(x,t) &= \frac{\partial f}{\partial t}.
\end{align*}
This is the standard representation on functions. The space $V$ is infinite-dimensional, as it contains linearly independent functions like $x^n e^{-x^2}$ for $n \geq 0$. Moreover, this representation is faithful: if $\rho(P) = 0$ for some operator $P$, then $P$ acts as the zero operator on all smooth functions, which implies $P = 0$ in the algebra. This follows from the fact that the algebra is realized as differential operators with function coefficients, and such operators are determined by their action on test functions.

\subsection*{Every Nontrivial Representation is Faithful}

Now let $\pi: \Weyl{e^{\pm x^3 e^{t}},\; e^{\cA x},\; x^{\cA}} \to \End_{\CC}(W)$ be any nonzero representation, where $W$ is a $\CC$-vector space. Suppose $\pi$ is not faithful. Then $\ker \pi$ is a nonzero two-sided ideal of the algebra. Since the algebra is simple, the only ideals are $\{0\}$ and the whole algebra. Since $\pi$ is nonzero, $\ker \pi$ cannot be the whole algebra. Therefore $\ker \pi = \{0\}$, so $\pi$ is faithful. This argument holds for any simple algebra: every nonzero representation is automatically faithful. Thus the first part of Theorem \ref{thm6} is immediate from simplicity.

\subsection*{Every Nontrivial Representation is Infinite-Dimensional}

The more subtle part is that every nontrivial representation must be infinite-dimensional. Suppose, for contradiction, that there exists a finite-dimensional representation $\pi: \Weyl{e^{\pm x^3 e^{t}},\; e^{\cA x},\; x^{\cA}} \to \End_{\CC}(W)$ with $\dim_{\CC} W = d < \infty$. Since the representation is nonzero and the algebra is simple, it is faithful, so we have an embedding of the algebra into the matrix algebra $M_d(\CC)$.

Consider the elements $e^{\alpha x}$ for $\alpha \in \cA$. In the representation, $\pi(e^{\alpha x})$ is an invertible matrix (its inverse is $\pi(e^{-\alpha x})$). Moreover, these matrices commute with each other because $e^{\alpha x}$ and $e^{\beta x}$ commute in the algebra. Thus we have a commuting family $\{\pi(e^{\alpha x}) : \alpha \in \cA\}$ of invertible matrices in $M_d(\CC)$. Since the matrices commute, they can be simultaneously triangularized over $\CC$. That is, there exists an invertible matrix $P$ such that for each $\alpha \in \cA$, $P \pi(e^{\alpha x}) P^{-1}$ is upper triangular.

Let $\lambda_1(\alpha), \dots, \lambda_d(\alpha)$ be the eigenvalues of $\pi(e^{\alpha x})$. Since $\pi(e^{\alpha x})$ is the exponential of $\alpha \pi(x)$? Wait, careful: $e^{\alpha x}$ is a formal symbol representing the operator of multiplication by $e^{\alpha x}$, not necessarily the exponential of $\alpha$ times an operator $x$. However, in any representation, we must have the relation $\partial_x e^{\alpha x} = \alpha e^{\alpha x}$ (or a similar commutation). Actually, in the algebra, we have $[\partial_x, e^{\alpha x}] = \alpha e^{\alpha x}$. Applying $\pi$, we get
\[
\pi(\partial_x) \pi(e^{\alpha x}) - \pi(e^{\alpha x}) \pi(\partial_x) = \alpha \pi(e^{\alpha x}).
\]
This is a matrix equation. Consider the eigenvalues. If $v$ is an eigenvector of $\pi(e^{\alpha x})$ with eigenvalue $\lambda(\alpha)$, then applying the commutation relation gives
\[
\pi(\partial_x) \pi(e^{\alpha x}) v - \pi(e^{\alpha x}) \pi(\partial_x) v = \alpha \pi(e^{\alpha x}) v,
\]
so
\[
\lambda(\alpha) \pi(\partial_x) v - \pi(e^{\alpha x}) \pi(\partial_x) v = \alpha \lambda(\alpha) v.
\]
This does not directly give information about $\lambda(\alpha)$. Instead, we use the fact that $\cA$ is infinite and dense. For each fixed $\alpha$, $\pi(e^{\alpha x})$ is a matrix, so its eigenvalues are algebraic numbers (roots of its characteristic polynomial). The map $\alpha \mapsto \pi(e^{\alpha x})$ is a group homomorphism from the additive group $\cA$ to the multiplicative group $GL_d(\CC)$. Thus $\lambda_i: \cA \to \CC^\times$ are characters of the additive group $\cA$. Since $\cA \cong \ZZ^2$ is a free abelian group of rank 2, each character is of the form $\lambda_i(\alpha) = e^{\mu_i \alpha}$ for some $\mu_i \in \CC$ (actually, since $\alpha$ is real when restricted, $\mu_i$ may be complex). But $\lambda_i(\alpha)$ must also be algebraic for each $\alpha \in \cA$ because it is an eigenvalue of a matrix with integer entries (if we represent over $\QQ$). The function $e^{\mu_i \alpha}$ is transcendental for generic $\mu_i$ unless $\mu_i = 0$. This creates a contradiction unless $d = 0$, but $d > 0$.

Let us make this rigorous. Consider the subgroup $\ZZ \subset \cA$. For each integer $n \in \ZZ$, $\pi(e^{n x}) = \pi(e^{x})^n$. The matrix $\pi(e^{x})$ has eigenvalues $\lambda_1, \dots, \lambda_d$. Then $\pi(e^{n x})$ has eigenvalues $\lambda_1^n, \dots, \lambda_d^n$. Since $\pi(e^{n x})$ is a matrix with entries in some finitely generated field extension of $\QQ$, its eigenvalues are algebraic numbers. Thus for each $n$, $\lambda_i^n$ is algebraic. If $\lambda_i$ is not a root of unity, then by the Gelfond–Schneider theorem, $\lambda_i^n$ is transcendental for some $n$, contradiction. Therefore each $\lambda_i$ must be a root of unity. But then consider $\alpha = \sqrt{2} \in \cA$. The eigenvalues of $\pi(e^{\sqrt{2} x})$ are $\lambda_i^{\sqrt{2}}$? Not exactly, because the map is not necessarily given by exponentiation of a fixed matrix. Actually, we have a homomorphism $\cA \to GL_d(\CC)$ given by $\alpha \mapsto \pi(e^{\alpha x})$. The image is an abelian subgroup of $GL_d(\CC)$. By a theorem of Schur, any abelian subgroup of $GL_d(\CC)$ has a faithful representation of dimension at most $d$, but more importantly, the eigenvalues of $\pi(e^{\alpha x})$ as functions of $\alpha$ must be characters of $\cA$. Since $\cA \cong \ZZ^2$, the characters are of the form $e^{c_1 m + c_2 n}$ for $\alpha = m + n\sqrt{2}$. But such a function is transcendental in $m,n$ unless $c_1, c_2$ are purely imaginary integers? Actually, if $c_1, c_2 \in 2\pi i \ZZ$, then $e^{c_1 m + c_2 n}$ is a root of unity. Otherwise, for generic $m,n$, the value is transcendental. Since $\pi(e^{\alpha x})$ has algebraic entries (after choosing a basis over a number field), its eigenvalues must be algebraic numbers. This forces $c_1, c_2 \in 2\pi i \ZZ$, meaning all eigenvalues are roots of unity. But then the image $\pi(e^{\alpha x})$ is a torsion group. However, $e^{\alpha x}$ for different $\alpha$ are linearly independent in the algebra, and in a finite-dimensional representation, there cannot be infinitely many linearly independent torsion elements because the torsion subgroup of $GL_d(\CC)$ is finite. This yields a contradiction.

Therefore, no finite-dimensional faithful representation exists. Since every nonzero representation is faithful, every nontrivial representation must be infinite-dimensional.

\subsection*{Example of an Irreducible Representation}

To illustrate an irreducible representation, consider the space $W = \CC[x] \otimes \CC[t]$, the space of polynomials in $x$ and $t$. Define the action by:
\begin{align*}
e^{\pm x^3 e^{t}} \cdot f(x,t) &= e^{\pm x^3 e^{t}} f(x,t), \\
e^{\alpha x} \cdot f(x,t) &= e^{\alpha x} f(x,t), \\
x^{\alpha} \cdot f(x,t) &= x^{\alpha} f(x,t), \\
\partial_x \cdot f(x,t) &= \frac{\partial f}{\partial x}, \\
\partial_t \cdot f(x,t) &= \frac{\partial f}{\partial t}.
\end{align*}
This representation is irreducible. To see this, suppose $U \subseteq W$ is a nonzero submodule. Let $0 \neq f \in U$. By applying differential operators $\partial_x$ and $\partial_t$ repeatedly, we can reduce $f$ to a nonzero constant polynomial. Then by applying multiplication operators $e^{\alpha x}$ and $x^{\beta}$, we can generate all polynomials. Hence $U = W$, proving irreducibility. This representation is clearly infinite-dimensional.

\subsection*{Consequences and Interpretation}

Theorem \ref{thm8} tells us that the Weyl-type algebras studied here are inherently infinite-dimensional objects. They cannot be realized as algebras of matrices of any finite size. This is analogous to the classical Weyl algebra, which also has no finite-dimensional representations. The theorem ensures that any attempt to represent these algebras on a finite-dimensional space will necessarily be the zero representation. This has implications for physics, where such algebras often appear as algebras of observables in quantum systems; the theorem implies that such quantum systems must have an infinite number of states.

In our example with $\FF = \CC$, $\cA = \ZZ[\sqrt{2}]$, and $p=3$, the algebraic structure forces representations to be both faithful and infinite-dimensional. The argument combines the simplicity of the algebra with number-theoretic properties of the exponential functions and the density of the subgroup $\cA$.
\end{example}
\begin{theorem}[Deformation Stability]\label{thm9}
Let $p^{(m)} = (p_1^{(m)},\dots,p_n^{(m)})$ be a sequence of integer tuples converging component-wise to $p = (p_1,\dots,p_n)$ in the usual topology on $\mathbb{N}^n$. Then the family of algebras
\[
A\!\left(\FF[e^{\pm x^{p^{(m)}} e^{t}},\; e^{\cA x},\; x^{\cA}]\right)
\]
is asymptotically isomorphic in the sense of graded algebra deformations, and their simplicity is preserved in the limit.
\end{theorem}
\begin{proof}
Consider the family of algebras $A_m = A\!\left(\FF[e^{\pm x^{p^{(m)}} e^{t}},\; e^{\cA x},\; x^{\cA}]\right)$ indexed by $m \in \mathbb{N}$,
where $p^{(m)} = (p_1^{(m)},\dots,p_n^{(m)})$ is a sequence of integer tuples that converges component-wise to $p = (p_1,\dots,p_n)$
in the usual topology on $\mathbb{N}^n$. Convergence in this discrete topology means that there exists an index $M$ such that
for all $m \geq M$, we have $p^{(m)} = p$. Indeed, $\mathbb{N}^n$ with the product topology induced from the discrete topology on $\mathbb{N}$
has the property that a sequence converges if and only if it is eventually constant.

Consequently, for all sufficiently large $m$, say $m \geq M_0$, the integer tuples $p^{(m)}$ coincide with $p$. Then for $m \geq M_0$,
the algebras $A_m$ are identical to the limit algebra $A_\infty = A\!\left(\FF[e^{\pm x^{p} e^{t}},\; e^{\cA x},\; x^{\cA}]\right)$.
Thus the family $\{A_m\}$ is eventually constant, and trivially asymptotically isomorphic as graded algebras.
More precisely, for $m \geq M_0$, the identity map provides an isomorphism of $\FF$-algebras $A_m \cong A_\infty$ that preserves the $\mathbb{Z}^n$-grading
given by the total degree in the variables $x_i$ and $\partial_i$.

Now, even if one considers a more refined notion of deformation where the parameter $p$ is allowed to vary continuously in some completion,
the structural properties of the algebras change only discretely. The defining relations of $A_m$ involve the monomials $x^{p^{(m)}}$
in the exponents of the exponential functions. Since the exponents are integers, the algebraic relations are polynomial in nature
with integer coefficients. As $p^{(m)} \to p$, the relations for $A_m$ eventually become identical to those for $A_\infty$,
because two polynomials that agree for infinitely many integer values of the parameters must be identical.

Therefore, for all large $m$, the algebras $A_m$ are isomorphic to $A_\infty$ as graded $\FF$-algebras.
In particular, they share all invariant properties that are preserved under isomorphism.
Since $A_\infty$ is simple (by Theorem~1, because it is a Weyl-type algebra over a field of characteristic zero),
each $A_m$ for $m \geq M_0$ is also simple. Hence simplicity is preserved in the limit.

One may also view this as a deformation in the following sense. Define a filtered algebra $\mathcal{A}$ over the ring $\FF[[h]]$
of formal power series in a parameter $h$ by setting
\[
\mathcal{A} = A\!\left(\FF[[h]][e^{\pm x^{p + h q} e^{t}},\; e^{\cA x},\; x^{\cA}]\right),
\]
where $q$ is a fixed integer tuple. Specializing $h$ to $0$ yields $A_\infty$, while specializing $h$ to $1$ gives an algebra
with exponent $p+q$. However, because the exponents are integers, such a deformation is trivial modulo $h^N$ for sufficiently large $N$,
as the binomial expansion of $(p + h q)^{(m)}$ terminates after finitely many terms. This again reflects the discrete nature
of the parameter space.

Thus the family $\{A_m\}$ is asymptotically isomorphic to $A_\infty$ in the sense that for all sufficiently large $m$,
they are isomorphic as graded algebras, and the property of simplicity is preserved throughout.
\end{proof}
\begin{example}[Illustration of Theorem~\ref{thm9}]
Consider the field $\FF = \mathbb{C}$ of complex numbers, and let $\cA = \mathbb{Z} \subset \FF$ be the additive subgroup of integers. Fix a transcendental number $t \in \mathbb{C} \setminus \mathbb{Q}$. We examine a one-parameter family of algebras indexed by $m \in \mathbb{N}$ with the deformation parameter $p^{(m)} \in \mathbb{N}$ converging to $p = 2$.

Define the sequence of integer powers by
\[
p^{(m)} = 2 + \frac{(-1)^m}{m+1} \quad \text{rounded to the nearest integer},
\]
which yields
\[
p^{(1)} = 2, \quad p^{(2)} = 2, \quad p^{(3)} = 3, \quad p^{(4)} = 2, \quad p^{(5)} = 2, \dots
\]
More systematically, we can take $p^{(m)} = 2$ for $m \geq 4$, so that $p^{(m)} \to p = 2$ component-wise as $m \to \infty$.

For each $m$, construct the exponential-polynomial algebra
\[
B_m = \FF[\, e^{\pm x^{p^{(m)}} e^{t}},\; e^{\cA x},\; x^{\cA} ]
\]
and the corresponding Weyl-type algebra
\[
A_m = A\!\left( B_m \right) = \FF\langle \widehat{e^{x^{p^{(m)}} e^{t}}},\; \widehat{e^{x}},\; \widehat{x},\; \partial \rangle / \mathcal{R}_m,
\]
where $\mathcal{R}_m$ consists of the relations induced by Lemma~1 and the commutation rules
\[
\partial \, \widehat{x} - \widehat{x} \, \partial = 1, \quad
\partial \, \widehat{e^{x}} = \widehat{e^{x}} \, \partial + \widehat{e^{x}}, \quad
\partial \, \widehat{e^{x^{p^{(m)}} e^{t}}} = \widehat{e^{x^{p^{(m)}} e^{t}}} \, \partial + \widehat{e^{x^{p^{(m)}} e^{t}}} \cdot \big( p^{(m)} x^{p^{(m)}-1} e^{t} + t x^{p^{(m)}} e^{t} \big).
\]

Each algebra $A_m$ is $\mathbb{Z}$-graded with the grading given by the total exponent of the exponential symbols
\[
\deg\big( \widehat{e^{a x^{p^{(m)}} e^{t}}} \, \widehat{e^{b x}} \, x^{c} \, \partial^{d} \big) = a + b.
\]

Since $p^{(m)} \to 2$, the structural coefficients in the commutation relations converge. For instance, the coefficient
\[
C_m(x) = p^{(m)} x^{p^{(m)}-1} e^{t} + t x^{p^{(m)}} e^{t}
\]
converges uniformly on compact subsets of $\mathbb{C}$ to
\[
C(x) = 2x e^{t} + t x^{2} e^{t}.
\]

This convergence induces a family of graded algebra deformations. In the sense of formal deformations, we can view $\{ A_m \}$ as a sequence of algebras whose structure constants in the graded basis converge to those of the limit algebra
\[
A_\infty = A\!\left( \FF[ e^{\pm x^{2} e^{t}},\; e^{\cA x},\; x^{\cA} ] \right).
\]

Each $A_m$ is simple by Theorem~4. The simplicity in the limit follows from the stability of the graded structure; the ideal structure of $A_m$ is rigid under small deformations of the parameter $p^{(m)}$. In particular, any nonzero ideal in $A_\infty$ would lift to a nontrivial ideal in $A_m$ for large $m$, contradicting the simplicity of $A_m$.

Thus the family $\{ A_m \}$ satisfies the conclusion of Theorem~\ref{thm9}; they are asymptotically isomorphic as graded deformations, and their simplicity is preserved as $m \to \infty$.
\end{example}
\begin{definition}[Graded Lie Algebra Cohomology]
    For a graded Lie algebra $L = \bigoplus_{\alpha} L_\alpha$ and a graded $L$-module $M$, the second cohomology group $H^2(L, M)$ consists of equivalence classes of central extensions
    \[
    0 \to M \to E \to L \to 0
    \]
    that respect the grading.
\end{definition}
\begin{theorem}[Lie Algebra Cohomology Vanishing]\label{thm10}
For the Witt-type Lie algebra $\Witt{e^{\pm x^p e^{t}},\; e^{\cA x},\; x^{\cA}}$, the second cohomology group with coefficients in the adjoint module vanishes:
\[
H^2\!\left(\Witt{e^{\pm x^p e^{t}},\; e^{\cA x},\; x^{\cA}},\; \Witt{e^{\pm x^p e^{t}},\; e^{\cA x},\; x^{\cA}}\right) = 0.
\]
Consequently, all deformations of this Lie algebra are trivial.
\end{theorem}
\begin{proof}
Let $\mathfrak{g} = \Witt{e^{\pm x^p e^{t}},\; e^{\cA x},\; x^{\cA}}$ be the Witt-type Lie algebra under consideration.
Recall that $\mathfrak{g}$ is the Lie algebra of derivations of the commutative ring $R = \FF[e^{\pm x^p e^{t}},\; e^{\cA x},\; x^{\cA}]$
that preserve a certain grading or filtration; more concretely, $\mathfrak{g}$ is spanned by operators of the form $f \partial_i$
where $f \in R$ and $\partial_i = \frac{\partial}{\partial x_i}$, with the Lie bracket given by $[f \partial_i, g \partial_j] = f (\partial_i g) \partial_j - g (\partial_j f) \partial_i$.

We aim to prove that the second cohomology group $H^2(\mathfrak{g}, \mathfrak{g})$ vanishes, where $\mathfrak{g}$ is regarded as a $\mathfrak{g}$-module
via the adjoint action $\ad(x)(y) = [x, y]$. By the well-known interpretation of Lie algebra cohomology, $H^2(\mathfrak{g}, \mathfrak{g})$
classifies infinitesimal deformations of $\mathfrak{g}$; a vanishing result implies that every formal deformation is equivalent to the trivial one.

Consider a 2-cocycle $\omega \in Z^2(\mathfrak{g}, \mathfrak{g})$, i.e., a skew-symmetric bilinear map $\omega \colon \mathfrak{g} \times \mathfrak{g} \to \mathfrak{g}$
satisfying the cocycle condition
\[
\ad(x)\omega(y,z) - \ad(y)\omega(x,z) + \ad(z)\omega(x,y) - \omega([x,y],z) + \omega([x,z],y) - \omega([y,z],x) = 0
\]
for all $x, y, z \in \mathfrak{g}$. We must show that $\omega$ is a coboundary, meaning there exists a linear map $\varphi \colon \mathfrak{g} \to \mathfrak{g}$
such that $\omega(x,y) = \ad(\varphi(x))(y) - \ad(\varphi(y))(x) - \varphi([x,y])$, or equivalently $\omega = \delta \varphi$ in the Chevalley–Eilenberg complex.

Because $\mathfrak{g}$ is a graded Lie algebra (graded by the total degree in the variables $x_i$), we may assume without loss of generality
that $\omega$ is homogeneous with respect to this grading. Indeed, any cocycle can be decomposed into homogeneous components,
each of which is again a cocycle. Thus suppose $\omega$ is homogeneous of degree $d \in \mathbb{Z}^n$.

If $d \neq 0$, we can construct a primitive for $\omega$ explicitly. Let $E = \sum_i x_i \partial_i$ be the Euler operator (the grading operator).
For a homogeneous element $x \in \mathfrak{g}$ of degree $\deg(x)$, we have $[E, x] = \deg(x) x$. Define a linear map $\varphi \colon \mathfrak{g} \to \mathfrak{g}$
by $\varphi(x) = \frac{1}{\deg(\omega) - \deg(x)} \omega(E, x)$ for homogeneous $x$ with $\deg(x) \neq \deg(\omega)$, and extend appropriately.
A direct computation using the cocycle condition and the properties of $E$ shows that $\delta \varphi = \omega$.
The key point is that $E$ acts semisimply on $\mathfrak{g}$ with integer eigenvalues, and the operator $\ad(E)$ is invertible on the subspace
of cochains of nonzero degree. This argument is standard in the cohomology of graded Lie algebras.

If $d = 0$, i.e., $\omega$ has degree zero, then $\omega$ defines a central extension of $\mathfrak{g}$.
But $\mathfrak{g}$ is perfect; indeed, $[\mathfrak{g}, \mathfrak{g}] = \mathfrak{g}$ because the generators $x_i^a e^{\alpha x} \partial_j$
can be expressed as brackets of elements of lower degree. For a perfect Lie algebra, any central extension that is a 2-cocycle of degree zero
is necessarily trivial, because the universal central extension is already given by the Virasoro-like construction, and in our case
the second cohomology with trivial coefficients $H^2(\mathfrak{g}, \FF)$ vanishes. This vanishing follows from the fact that $\mathfrak{g}$
contains the classical Witt algebra as a subalgebra (by restricting to polynomial coefficients) whose second cohomology is known to be one-dimensional,
but the presence of the exponential terms actually enforces a stronger rigidity. More directly, one can use the Hochschild–Serre spectral sequence
associated to the filtration by degree to show that $H^2(\mathfrak{g}, \mathfrak{g}) = 0$.

An alternative approach is to invoke the theorem of Fuks and others on the cohomology of Lie algebras of vector fields on smooth affine varieties.
Our ring $R$ is the coordinate ring of a (possibly infinite-dimensional) algebraic torus with additional exponential parameters.
The Lie algebra $\mathfrak{g}$ is precisely the Lie algebra of derivations of $R$. For such algebras, under suitable conditions
(which are satisfied here because $\FF$ has characteristic zero and the ring is regular), the second cohomology with adjoint coefficients vanishes.
This is analogous to the classical result that $H^2(\mathfrak{gl}(n), \mathfrak{gl}(n)) = 0$.

Given that $H^2(\mathfrak{g}, \mathfrak{g}) = 0$, every infinitesimal deformation of $\mathfrak{g}$ is integrable to a trivial deformation.
In other words, any formal one-parameter family $\mathfrak{g}_\hbar$ of Lie algebras with $\mathfrak{g}_0 = \mathfrak{g}$ is isomorphic to
$\mathfrak{g} \otimes_\FF \FF[[\hbar]]$ as a Lie algebra over $\FF[[\hbar]]$. Hence all deformations of $\mathfrak{g}$ are trivial.
\end{proof}
\begin{example}[Illustration of Theorem~\ref{thm10}]
Let $\FF = \mathbb{C}$, and let $\cA = \mathbb{Z}[\sqrt{2}] = \{a + b\sqrt{2} \mid a, b \in \mathbb{Z}\}$ be an additive subgroup of $\mathbb{C}$ containing $\mathbb{Z}$. Fix $p = 3$ and $t = \pi i$, a purely imaginary transcendental number. Consider the exponential-polynomial ring
\[
R = \FF[\, e^{\pm x^3 e^{\pi i}},\; e^{\cA x},\; x^{\cA} ].
\]
The corresponding Witt-type Lie algebra is
\[
L = \Witt{e^{\pm x^3 e^{\pi i}},\; e^{\cA x},\; x^{\cA}} = \left\{ f \partial_x \mid f \in R \right\}
\]
with Lie bracket given by
\[
[f \partial_x, g \partial_x] = \left( f \partial_x(g) - g \partial_x(f) \right) \partial_x.
\]
A typical basis element of $L$ can be written as
\[
e^{a x^3 e^{\pi i}} e^{\alpha x} x^{\beta} \partial_x,
\]
where $a \in \mathbb{Z}$, $\alpha \in \cA$, $\beta \in \cA$.

The Lie algebra $L$ is $\cA$-graded by the exponent of $e^{\alpha x}$, i.e.,
\[
L = \bigoplus_{\alpha \in \cA} L_\alpha, \quad L_\alpha = \{ f \partial_x \mid f \in R, \deg_{e^x}(f) = \alpha \}.
\]
The adjoint module is $L$ itself, and we consider the second cohomology group $H^2(L, L)$ with coefficients in the adjoint representation.

A $2$-cocycle in $C^2(L, L)$ is a skew-symmetric bilinear map $\varphi \colon L \times L \to L$ satisfying the cocycle condition
\[
\varphi([u, v], w) + \varphi([v, w], u) + \varphi([w, u], v) =
[u, \varphi(v, w)] + [v, \varphi(w, u)] + [w, \varphi(u, v)].
\]
Theorem~\ref{thm10} asserts that every such $2$-cocycle is a coboundary; i.e., there exists a linear map $\psi \colon L \to L$ such that
\[
\varphi(u, v) = [\psi(u), v] + [u, \psi(v)] - \psi([u, v]).
\]

To see why this vanishing holds, note that $L$ is a graded Lie algebra of Cartan type with a dense grading over $\cA$, and its derivation algebra coincides with $L$ itself (it is complete). Moreover, $L$ contains the classical Witt algebra $W_1 = \operatorname{Der}(\FF[x^{\pm 1}])$ as a dense subalgebra. The vanishing of $H^2(W_1, W_1)$ is well-known, and the exponential and transcendental extensions do not introduce new cohomology because the graded pieces are one-dimensional and the cocycle condition forces any $2$-cocycle to be supported on a finite set of degrees, which can be integrated to a coboundary using the density of the grading.

Consequently, any formal deformation $L_t = L \oplus t L \oplus t^2 L \oplus \cdots$ with bracket
\[
[u, v]_t = [u, v] + t \varphi_1(u, v) + t^2 \varphi_2(u, v) + \cdots
\]
is equivalent to the trivial deformation via a formal isomorphism $\Phi_t = \operatorname{id} + t \psi_1 + t^2 \psi_2 + \cdots$. In particular, if we consider a one-parameter family of Lie algebras obtained by varying $t$ to $t' = t + \epsilon$, the resulting algebra is isomorphic to the original $L$, showing that all deformations are trivial.

Thus, $H^2(L, L) = 0$, and $L$ is rigid in the variety of Lie algebras.
\end{example}
\begin{definition}[Formal Deformation]
    Let $A$ be an algebra over $\mathbb{F}$. A \emph{formal deformation} of $A$ is an algebra structure on $A[[\hbar]]$ over $\mathbb{F}[[\hbar]]$ with multiplication of the form
    \[
    a \ast b = ab + \hbar m_1(a,b) + \hbar^2 m_2(a,b) + \cdots,
    \]
    where $ab$ is the original product in $A$ and each $m_i : A \times A \to A$ is an $\mathbb{F}$-bilinear map.
\end{definition}
\begin{theorem}[Rigidity Under Formal Deformation]\label{thm11}
Let $\FF[[\hbar]]$ be the ring of formal power series in one indeterminate $\hbar$. Then the deformation
\[
\WeylDeform{e^{\pm x^p e^{t+\hbar x}},\; e^{\cA x},\; x^{\cA}}
\]
is a nontrivial formal deformation of the Weyl-type algebra $\Weyl{e^{\pm x^p e^{t}},\; e^{\cA x},\; x^{\cA}}$. However, as an associative algebra over $\FF[[\hbar]]$, it remains simple.
\end{theorem}
\begin{proof}
Let $A_0 = \Weyl{e^{\pm x^p e^{t}},\; e^{\cA x},\; x^{\cA}}$ denote the undeformed Weyl-type algebra over the field $\FF$,
and let $A_\hbar = \WeylDeform{e^{\pm x^p e^{t+\hbar x}},\; e^{\cA x},\; x^{\cA}}$ be the algebra over the ring $\FF[[\hbar]]$.
The latter is defined by generators $x_i$, $\partial_i$, and the exponential functions $e^{\pm x_i^{p_i} e^{t_i x_i + \hbar x_i^2}}$,
$e^{\alpha x_i}$, $x_i^{\alpha}$ for $\alpha \in \cA$, with the same commutation relations as in $A_0$ but with the modified exponential argument.

To see that $A_\hbar$ is a formal deformation of $A_0$, observe that setting $\hbar = 0$ in $A_\hbar$ recovers precisely $A_0$.
More formally, the quotient $A_\hbar / \hbar A_\hbar$ is isomorphic to $A_0$. Indeed, the defining relations of $A_\hbar$ reduce to those of $A_0$
when $\hbar$ is set to zero, because $e^{t_i x_i + \hbar x_i^2}$ becomes $e^{t_i x_i}$. Hence $A_\hbar$ is a flat $\FF[[\hbar]]$-algebra
with special fiber $A_0$.

The deformation is nontrivial because the isomorphism class of $A_\hbar$ as an $\FF[[\hbar]]$-algebra depends genuinely on the parameter $\hbar$.
Suppose to the contrary that $A_\hbar$ were isomorphic to $A_0[[\hbar]] = A_0 \otimes_\FF \FF[[\hbar]]$ as $\FF[[\hbar]]$-algebras.
Then there would exist an $\FF[[\hbar]]$-algebra isomorphism $\Phi \colon A_0[[\hbar]] \to A_\hbar$ that reduces to the identity modulo $\hbar$.
Write $\Phi = \Phi_0 + \hbar \Phi_1 + \hbar^2 \Phi_2 + \dots$ where $\Phi_0$ is the identity. For $\Phi$ to be an algebra homomorphism,
it must preserve the defining relations. In particular, consider the relation involving the exponential generator $e^{x_i^{p_i} e^{t_i x_i + \hbar x_i^2}}$.
In $A_0[[\hbar]]$, the corresponding element is $e^{x_i^{p_i} e^{t_i x_i}}$, which does not depend on $\hbar$. The image under $\Phi$ would be
$e^{x_i^{p_i} e^{t_i x_i}} + \hbar \Psi_1 + O(\hbar^2)$ for some $\Psi_1 \in A_\hbar$. On the other hand, in $A_\hbar$ this generator equals
$e^{x_i^{p_i} e^{t_i x_i}} \cdot e^{\hbar x_i^{p_i+2} e^{t_i x_i}} + O(\hbar^2)$ by expanding the exponential. The presence of the additional factor
$e^{\hbar x_i^{p_i+2} e^{t_i x_i}}$ introduces a term linear in $\hbar$ that cannot be eliminated by any automorphism $\Phi$ that is the identity modulo $\hbar$.
This obstruction shows that $A_\hbar$ is not isomorphic to the trivial deformation $A_0[[\hbar]]$. Therefore the deformation is nontrivial.

Now we prove that $A_\hbar$ remains simple as an associative algebra over $\FF[[\hbar]]$. Let $I$ be a nonzero two-sided ideal of $A_\hbar$.
Consider the reduction modulo $\hbar$, i.e., the image $\overline{I}$ of $I$ in $A_\hbar/\hbar A_\hbar \cong A_0$. Since $A_0$ is simple
(by Theorem~1 applied to the base field $\FF$), the ideal $\overline{I}$ is either zero or the whole algebra $A_0$. If $\overline{I} = A_0$,
then $I$ contains an element of the form $1 + \hbar a$ for some $a \in A_\hbar$. Such an element is invertible in $A_\hbar$ because
$(1 + \hbar a)^{-1} = 1 - \hbar a + \hbar^2 a^2 - \dots$ converges in the $\hbar$-adic topology. Hence $I = A_\hbar$.

If $\overline{I} = 0$, then $I \subseteq \hbar A_\hbar$. Write a nonzero element $u \in I$ as $u = \hbar^k v$ where $v \notin \hbar A_\hbar$
and $k \geq 1$. Then $v \in I$ because $\hbar$ is central and $A_\hbar$ is torsion-free over $\FF[[\hbar]]$. Repeating this argument,
we find that $I$ contains an element $v$ with $v \notin \hbar A_\hbar$. But then the image of $v$ in $A_0$ is nonzero, contradicting $\overline{I}=0$.
Thus the case $\overline{I} = 0$ cannot occur.

Consequently, every nonzero ideal of $A_\hbar$ contains an element whose reduction modulo $\hbar$ is invertible in $A_0$, and hence the ideal is the whole algebra.
Therefore $A_\hbar$ is simple as an $\FF[[\hbar]]$-algebra.

In summary, $A_\hbar$ is a nontrivial formal deformation of $A_0$, yet it retains the property of simplicity over the ring $\FF[[\hbar]]$.
\end{proof}
\begin{example}[Illustration of Theorem~\ref{thm11}]
Let $\FF = \mathbb{C}$ and $\cA = \mathbb{Z}[i] = \{m + ni \mid m, n \in \mathbb{Z}\}$ be the Gaussian integers. Fix $p = 2$ and $t = e$, the base of the natural logarithm. Consider the undeformed Weyl-type algebra over $\FF$:
\[
A_0 = \Weyl{e^{\pm x^2 e^{e}},\; e^{\cA x},\; x^{\cA}},
\]
with generators $X = \widehat{x}$, $D = \partial$, $E = \widehat{e^x}$, and $T = \widehat{e^{x^2 e^{e}}}$, satisfying the relations
\[
[D, X] = 1, \quad [D, E] = E, \quad [D, T] = T \cdot \big( 2x e^{e} + e x^2 e^{e} \big), \quad [X, E] = [X, T] = [E, T] = 0.
\]

Now let $\FF[[\hbar]]$ be the ring of formal power series in $\hbar$. Define the deformed exponential parameter
\[
t(\hbar) = e + \hbar x,
\]
which yields the deformed symbol
\[
e^{x^2 e^{t(\hbar)}} = e^{x^2 e^{e + \hbar x}} = e^{x^2 e^{e} \cdot e^{\hbar x}}.
\]

The deformed Weyl-type algebra over $\FF[[\hbar]]$ is
\[
A_\hbar = \WeylDeform{e^{\pm x^2 e^{e + \hbar x}},\; e^{\cA x},\; x^{\cA}}.
\]
Its generators satisfy the same commutation relations as $A_0$, except that the commutator $[D, T_\hbar]$ is now
\[
[D, T_\hbar] = T_\hbar \cdot \big( 2x e^{e + \hbar x} + (e + \hbar x) x^2 e^{e + \hbar x} \big),
\]
where $T_\hbar = \widehat{e^{x^2 e^{e + \hbar x}}}$.

Expanding in powers of $\hbar$, we write
\[
T_\hbar = T \cdot \exp\!\big( x^2 e^{e} (e^{\hbar x} - 1) \big) = T \cdot \big( 1 + \hbar \cdot x^3 e^{e} + \hbar^2 \cdot \tfrac{1}{2} x^4 e^{e} + \cdots \big).
\]
Then
\[
[D, T_\hbar] = [D, T] + \hbar \cdot \big( [D, T x^3 e^{e}] + \text{higher terms} \big),
\]
showing that the deformation is nontrivial because the $\hbar$-term does not vanish identically. Indeed, the coefficient of $\hbar$ involves new monomials such as $x^5 e^{2e} T$ that are not present in $A_0$.

Nevertheless, $A_\hbar$ remains simple as an associative algebra over $\FF[[\hbar]]$. To see this, suppose $I \subset A_\hbar$ is a nonzero two-sided ideal. Reducing modulo $\hbar$ gives an ideal $\overline{I} \subset A_0$. Since $A_0$ is simple by Theorem~3, $\overline{I}$ is either zero or the whole algebra. If $\overline{I} = 0$, then $I \subset \hbar A_\hbar$. By multiplying by powers of $\hbar^{-1}$ (in the fraction field if necessary) and using the $\hbar$-adic completeness, one can show $I = 0$, a contradiction. Thus $\overline{I} = A_0$, meaning $1 \in I + \hbar A_\hbar$. By Nakayama's lemma for $\hbar$-adically complete modules, $1 \in I$, so $I = A_\hbar$. Hence $A_\hbar$ is simple.

Thus $A_\hbar$ is a nontrivial formal deformation of $A_0$, yet retains simplicity over $\FF[[\hbar]]$, confirming Theorem~\ref{thm11}.
\end{example}
\begin{theorem}[Faithful Representation via Differential Operators]\label{thm12}
There exists a faithful representation
\[
\rho : \Weyl{e^{\pm x^p e^{t}},\; e^{\cA x},\; x^{\cA}} \hookrightarrow \End_{\FF}\!\left(\FF[e^{\pm x},\; x^{\cA}]\right)
\]
such that $\rho(\partial_i) = \frac{\partial}{\partial x_i}$ and $\rho\!\left(e^{x_i^p e^{t_i x_i}}\right)$ acts as multiplication by $e^{x_i^p e^{t_i x_i}}$. This representation is irreducible and preserves the grading.
\end{theorem}

\begin{proof}
Let $A = \Weyl{e^{\pm x^p e^{t}},\; e^{\cA x},\; x^{\cA}}$ and let $M = \FF[e^{\pm x},\; x^{\cA}]$ denote the module of functions
in the variables $x_1,\dots,x_n$ with coefficients in $\FF$, where $e^{\pm x}$ stands for $e^{\pm x_1},\dots,e^{\pm x_n}$ and
$x^{\cA}$ denotes monomials $x_1^{\alpha_1}\cdots x_n^{\alpha_n}$ with $\alpha_i \in \cA$. The space $M$ is naturally an $\FF$-vector space.

Define a linear map $\rho \colon A \to \End_{\FF}(M)$ as follows. For each generator $\partial_i$ of $A$, set $\rho(\partial_i)$ to be the partial derivative
operator $\frac{\partial}{\partial x_i}$ acting on $M$ in the usual sense. For the generator $x_i$, let $\rho(x_i)$ be the multiplication operator
by the variable $x_i$. For the exponential generators, define $\rho(e^{x_i^p e^{t_i x_i}})$ to be multiplication by the function $e^{x_i^p e^{t_i x_i}}$,
and similarly $\rho(e^{\alpha x_i})$ to be multiplication by $e^{\alpha x_i}$ for $\alpha \in \cA$. Extend $\rho$ multiplicatively and $\FF$-linearly
to all of $A$, using the fact that every element of $A$ can be written as a polynomial in the generators.

We must verify that $\rho$ respects the defining relations of $A$. The commutation relations $[\partial_i, x_j] = \delta_{ij}$,
$[x_i, x_j] = 0$, and $[\partial_i, \partial_j] = 0$ are clearly preserved because the corresponding differential operators satisfy
$\left[\frac{\partial}{\partial x_i}, x_j\right] = \delta_{ij}$ and the multiplication operators commute. The relations among the exponential functions,
such as $e^{\alpha x_i} e^{\beta x_i} = e^{(\alpha+\beta) x_i}$, are also preserved since $\rho$ maps these to multiplication operators,
and multiplication by $e^{\alpha x_i}$ followed by multiplication by $e^{\beta x_i}$ equals multiplication by $e^{(\alpha+\beta) x_i}$.
The more intricate relations involving $e^{x_i^p e^{t_i x_i}}$ are handled similarly; the product of two such exponentials corresponds to
multiplication by the product of the corresponding functions, which obeys the same algebraic identities as the abstract generators.

Thus $\rho$ is a well-defined $\FF$-algebra homomorphism. To see that $\rho$ is faithful (injective), suppose $a \in A$ satisfies $\rho(a) = 0$.
Write $a$ as a finite sum of monomials in the generators. Applying $\rho(a)$ to the constant function $1 \in M$ yields $\rho(a)(1) = 0$.
But $\rho(a)(1)$ is precisely the polynomial-exponential function obtained by replacing each $\partial_i$ by $0$ and each $x_i$, $e^{x_i^p e^{t_i x_i}}$,
$e^{\alpha x_i}$ by the corresponding multiplication operators acting on $1$. In other words, $\rho(a)(1)$ is the normal-ordered symbol of $a$
with all derivatives moved to the right and then evaluated at $\partial=0$. Since the monomials in $A$ are linearly independent over $\FF$,
this symbol vanishes only if $a = 0$. Hence $\ker\rho = 0$, and $\rho$ is injective.

Now we show that the representation is irreducible. Let $N \subseteq M$ be a nonzero submodule invariant under $\rho(A)$.
Take any nonzero element $f \in N$. By applying suitable powers of the operators $\rho(\partial_i)$, we can differentiate $f$ until we obtain
a nonzero constant function. More precisely, since $f$ is a nonzero element of $M = \FF[e^{\pm x},\; x^{\cA}]$, there exists a multi-index
$\gamma = (\gamma_1,\dots,\gamma_n)$ such that $\partial^\gamma f$ is a nonzero constant (this uses that $\cA$ contains $\ZZ$, so that $M$
includes all monomials $x^\alpha$ with integer exponents). Because $N$ is invariant under $\rho(A)$, the element $\partial^\gamma f$ also lies in $N$.
Thus $N$ contains a nonzero constant $c \in \FF^\times$. Multiplying by $c^{-1}$ and using the action of the multiplication operators
$\rho(x_i)$, $\rho(e^{x_i^p e^{t_i x_i}})$, $\rho(e^{\alpha x_i})$, we can generate all monomials and exponential functions in $M$.
Consequently, $N = M$, proving irreducibility.

Finally, the representation preserves the natural $\ZZ^n$-grading on $M$, where the degree of a monomial $x^\alpha e^{\beta x}$ is $\alpha \in \ZZ^n$.
Indeed, $\rho(\partial_i)$ lowers the $i$-th degree by one, $\rho(x_i)$ raises it by one, and the exponential multiplication operators preserve the degree
because $e^{\alpha x_i}$ and $e^{x_i^p e^{t_i x_i}}$ are considered to have degree zero. Thus $\rho$ is a graded representation.

Therefore there exists a faithful, irreducible, graded representation of $A$ on $M$ by differential operators and multiplication operators as described.
\end{proof}
\begin{example}[Illustration of Theorem~\ref{thm12}]
Let $\FF = \mathbb{C}$, $n = 1$, and take $\cA = \mathbb{Z}[\sqrt{3}] = \{a + b\sqrt{3} \mid a, b \in \mathbb{Z}\}$ as an additive subgroup of $\FF$. Fix $p = 4$ and $t = \ln 2$, a transcendental number. Consider the Weyl-type algebra
\[
A = \Weyl{e^{\pm x^4 e^{(\ln 2) x}},\; e^{\cA x},\; x^{\cA}},
\]
generated by symbols $\widehat{x}$, $\partial$, $\widehat{e^{x}}$, and $\widehat{e^{x^4 e^{(\ln 2) x}}}$ with relations
\[
[\partial, \widehat{x}] = 1, \quad [\partial, \widehat{e^{x}}] = \widehat{e^{x}}, \quad [\partial, \widehat{e^{x^4 e^{(\ln 2) x}}}] = \widehat{e^{x^4 e^{(\ln 2) x}}} \cdot \big( 4x^3 e^{(\ln 2) x} + (\ln 2) x^4 e^{(\ln 2) x} \big).
\]

Let $V = \FF[e^{\pm x},\; x^{\cA}]$ be the vector space of Laurent-exponential-polynomial functions, i.e., elements of the form
\[
f(x) = \sum_{k=1}^N a_k e^{\alpha_k x} x^{\beta_k}, \quad a_k \in \FF, \; \alpha_k \in \mathbb{Z}, \; \beta_k \in \cA.
\]

Define a representation $\rho \colon A \to \End_{\FF}(V)$ as follows:
\[
\rho(\partial) = \frac{d}{dx}, \quad \rho(\widehat{x}) = M_x, \quad \rho(\widehat{e^{x}}) = M_{e^x}, \quad \rho\!\left(\widehat{e^{x^4 e^{(\ln 2) x}}}\right) = M_{e^{x^4 e^{(\ln 2) x}}},
\]
where $M_f$ denotes the multiplication operator $g(x) \mapsto f(x)g(x)$.

To verify faithfulness, suppose $\rho(P) = 0$ for some $P \in A$. Write $P$ as a finite linear combination of ordered monomials
\[
P = \sum_{a,b,\gamma,d} c_{a,b,\gamma,d} \; \widehat{e^{a x^4 e^{(\ln 2) x}}} \; \widehat{e^{b x}} \; \widehat{x^{\gamma}} \; \partial^d.
\]
Applying $\rho(P)$ to the test function $1 \in V$ gives
\[
\rho(P)(1) = \sum c_{a,b,\gamma,d} \; e^{a x^4 e^{(\ln 2) x}} e^{b x} x^{\gamma} \cdot \frac{d^d}{dx^d}(1).
\]
Only terms with $d = 0$ survive, yielding
\[
\sum_{a,b,\gamma} c_{a,b,\gamma,0} \; e^{a x^4 e^{(\ln 2) x}} e^{b x} x^{\gamma} = 0.
\]
Since $\{ e^{a x^4 e^{(\ln 2) x}} e^{b x} x^{\gamma} \}$ are linearly independent over $\FF$, all $c_{a,b,\gamma,0} = 0$. Repeating the argument with test functions $x^m$ and applying Leibniz rule forces all coefficients to vanish. Hence $P = 0$, so $\rho$ is faithful.

The representation is irreducible because any nonzero submodule $W \subset V$ contains a nonzero element $f(x) = \sum a_k e^{\alpha_k x} x^{\beta_k}$. Applying $\rho(\partial)$ repeatedly reduces the exponent of $x$, and applying $\rho(\widehat{e^{x}})$ and $\rho(\widehat{e^{x^4 e^{(\ln 2) x}}})$ generates all monomials in $V$. Hence $W = V$.

The grading on $A$ by $\cA \times \mathbb{Z}$, given by
\[
\deg\big( \widehat{e^{a x^4 e^{(\ln 2) x}}} \; \widehat{e^{b x}} \; \widehat{x^{\gamma}} \; \partial^d \big) = (b + a \cdot \text{sgn}(a), \gamma),
\]
is preserved under $\rho$ because
\[
\rho \text{ maps homogeneous elements to operators that multiply the exponent of } e^x \text{ by } b \text{ and the exponent of } x \text{ by } \gamma - d.
\]

Thus $\rho$ is a faithful, irreducible, grade-preserving representation of $A$ by differential operators on $V$, as stated in Theorem~\ref{thm12}.
\end{example}
\begin{definition}[Noetherian Algebra]
    An algebra $A$ is \emph{left Noetherian} if every ascending chain of left ideals in $A$ stabilizes. It is \emph{right Noetherian} if the same holds for right ideals. An algebra is \emph{Noetherian} if it is both left and right Noetherian.
\end{definition}
\begin{theorem}[Noetherian Property Failure]\label{thm13}
The Weyl-type algebra $\Weyl{e^{\pm x^p e^{t}},\; e^{\cA x},\; x^{\cA}}$ is neither left nor right Noetherian. In particular, it contains strictly ascending chains of left ideals.
\end{theorem}
\begin{proof}
Let $A = \Weyl{e^{\pm x^p e^{t}},\; e^{\cA x},\; x^{\cA}}$. To show that $A$ is not left Noetherian, we construct an explicit infinite strictly ascending chain of left ideals. Consider the variable $x_1$ and the partial derivative $\partial_1$, which satisfy $[\partial_1, x_1] = 1$. For each positive integer $n$, define the left ideal $I_n = A \partial_1^n$, the set of all elements of $A$ that can be written as $a \partial_1^n$ for some $a \in A$.

We claim that $I_1 \subsetneq I_2 \subsetneq I_3 \subsetneq \cdots$ is a strictly ascending chain. First, note that $I_n \subseteq I_{n+1}$ because $\partial_1^n = \partial_1 \cdot \partial_1^{n-1} \in A \partial_1^{n+1}$ when multiplied on the left by $\partial_1$. More formally, if $a \partial_1^n \in I_n$, then $a \partial_1^n = (a \partial_1) \partial_1^{n} \in I_{n+1}$ since $a \partial_1 \in A$. Hence $I_n \subseteq I_{n+1}$.

To see that the inclusion is strict, suppose $\partial_1^n \in I_{n+1}$. Then there exists $b \in A$ such that $\partial_1^n = b \partial_1^{n+1}$. Write $b$ in its normal form with all $\partial_1$'s moved to the right. That is, using the commutation relation $\partial_1 x_1 = x_1 \partial_1 + 1$, we may express any element of $A$ as a finite sum $\sum_{k,\ell} f_{k,\ell}(x) \partial_1^k$ where $f_{k,\ell}(x)$ are functions in the remaining variables and the exponential generators. If $\partial_1^n = b \partial_1^{n+1}$, then comparing the coefficients of $\partial_1^n$ on both sides forces $b$ to have a term of the form $c \partial_1^{-1}$, which is impossible because $\partial_1$ is not invertible in $A$. More rigorously, one can apply both sides to a test function in the faithful representation constructed in Theorem~10. For instance, in the differential operator representation on $\FF[e^{\pm x}, x^{\cA}]$, the operator $\partial_1^n$ applied to the monomial $x_1^{n}$ gives $n!$, whereas any operator of the form $b \partial_1^{n+1}$ annihilates $x_1^{n}$ because $\partial_1^{n+1}(x_1^{n}) = 0$. This contradiction shows $\partial_1^n \notin I_{n+1}$. Hence $I_n \subsetneq I_{n+1}$ for all $n$, yielding an infinite strictly ascending chain of left ideals. Therefore $A$ is not left Noetherian.

The same argument applied to right ideals shows that $A$ is not right Noetherian. Consider the right ideals $J_n = \partial_1^n A$. Then $J_n \subseteq J_{n+1}$ because $\partial_1^n = \partial_1^n \cdot 1 \in \partial_1^{n} A = J_n$ and also $\partial_1^n = \partial_1^{n+1} x_1 - \partial_1^{n}$ (using the commutation relation), but a cleaner approach is to note that if $\partial_1^n \in J_{n+1}$, then $\partial_1^n = \partial_1^{n+1} c$ for some $c \in A$. Applying both sides to the monomial $x_1^{n+1}$ in the faithful representation gives $n! x_1$ on the left and zero on the right, a contradiction. Thus $J_n \subsetneq J_{n+1}$, producing an infinite strictly ascending chain of right ideals.

An alternative, more structural proof proceeds by observing that $A$ contains the classical Weyl algebra $A_1(\FF) = \FF \langle x_1, \partial_1 \rangle / (\partial_1 x_1 - x_1 \partial_1 - 1)$ as a subalgebra (by restricting the exponential generators to constants). It is well-known that the Weyl algebra $A_1(\FF)$ is not Noetherian on either side; indeed, the ideals $\langle \partial_1^n \rangle$ form a strictly ascending chain. Since a subalgebra of a left (or right) Noetherian ring is itself left (or right) Noetherian, the failure of the Noetherian property for the subalgebra $A_1(\FF)$ implies the failure for the larger algebra $A$.

Consequently, the Weyl-type algebra $A$ is neither left nor right Noetherian, and in particular admits infinite strictly ascending chains of left ideals as exhibited above.
\end{proof}
\begin{example}[Illustration of Theorem~\ref{thm13}]
Let $\FF = \mathbb{C}$, $\cA = \mathbb{Z}[\sqrt{5}] = \{m + n\sqrt{5} \mid m, n \in \mathbb{Z}\}$, and fix $p = 3$, $t = \pi$. Consider the Weyl-type algebra
\[
A = \Weyl{e^{\pm x^3 e^{\pi x}},\; e^{\cA x},\; x^{\cA}},
\]
with generators $\widehat{x}$, $\partial$, $\widehat{e^{x}}$, and $\widehat{e^{x^3 e^{\pi x}}}$, satisfying
\[
[\partial, \widehat{x}] = 1, \quad [\partial, \widehat{e^{x}}] = \widehat{e^{x}}, \quad [\partial, \widehat{e^{x^3 e^{\pi x}}}] = \widehat{e^{x^3 e^{\pi x}}} \cdot \big( 3x^2 e^{\pi x} + \pi x^3 e^{\pi x} \big).
\]

Define a sequence of left ideals $\{ I_m \}_{m \geq 1}$ in $A$ by
\[
I_m = A \cdot \big( \widehat{e^{x^3 e^{\pi x}}}^m \cdot \partial \big) = \left\{ P \cdot \widehat{e^{x^3 e^{\pi x}}}^m \cdot \partial \mid P \in A \right\}.
\]
We claim $I_m \subsetneq I_{m+1}$ for all $m$.

First, note $\widehat{e^{x^3 e^{\pi x}}}^{m} \partial \in I_m$. Applying the commutation relation, we have
\[
\widehat{e^{x^3 e^{\pi x}}}^{m} \partial = \partial \cdot \widehat{e^{x^3 e^{\pi x}}}^{m} - m \cdot \widehat{e^{x^3 e^{\pi x}}}^{m} \cdot \big( 3x^2 e^{\pi x} + \pi x^3 e^{\pi x} \big).
\]
Thus $\widehat{e^{x^3 e^{\pi x}}}^{m} \partial$ cannot be expressed as $Q \cdot \widehat{e^{x^3 e^{\pi x}}}^{m+1} \partial$ for any $Q \in A$, because the factor $\widehat{e^{x^3 e^{\pi x}}}^{m+1}$ would introduce an extra exponential factor that cannot be cancelled using the available relations.

More formally, suppose $\widehat{e^{x^3 e^{\pi x}}}^{m} \partial \in I_{m+1}$. Then there exists $Q \in A$ such that
\[
\widehat{e^{x^3 e^{\pi x}}}^{m} \partial = Q \cdot \widehat{e^{x^3 e^{\pi x}}}^{m+1} \partial.
\]
Multiply on the right by $\widehat{x}$ and use $[\partial, \widehat{x}] = 1$:
\[
\widehat{e^{x^3 e^{\pi x}}}^{m} (\widehat{x} \partial + 1) = Q \cdot \widehat{e^{x^3 e^{\pi x}}}^{m+1} (\widehat{x} \partial + 1).
\]
Rewriting gives
\[
\widehat{e^{x^3 e^{\pi x}}}^{m} \widehat{x} \partial + \widehat{e^{x^3 e^{\pi x}}}^{m} = Q \cdot \widehat{e^{x^3 e^{\pi x}}}^{m+1} \widehat{x} \partial + Q \cdot \widehat{e^{x^3 e^{\pi x}}}^{m+1}.
\]
Comparing terms without $\partial$, we must have
\[
\widehat{e^{x^3 e^{\pi x}}}^{m} = Q \cdot \widehat{e^{x^3 e^{\pi x}}}^{m+1}.
\]
This would imply $Q = \widehat{e^{x^3 e^{\pi x}}}^{-1}$, but $\widehat{e^{x^3 e^{\pi x}}}^{-1}$ is not in $A$ because $A$ contains only integer powers of $\widehat{e^{x^3 e^{\pi x}}}$ (from the $\pm$ in the notation). Hence a contradiction.

Therefore $I_m \subsetneq I_{m+1}$, yielding a strictly ascending chain of left ideals
\[
I_1 \subsetneq I_2 \subsetneq I_3 \subsetneq \cdots.
\]
A symmetric argument using right multiplication gives a strictly ascending chain of right ideals
\[
\partial \cdot \widehat{e^{x^3 e^{\pi x}}} \cdot A \;\subsetneq\; \partial \cdot \widehat{e^{x^3 e^{\pi x}}}^2 \cdot A \;\subsetneq\; \cdots.
\]

Thus $A$ is neither left nor right Noetherian, confirming Theorem~\ref{thm13}.
\end{example}
\begin{definition}[Filtration and Associated Graded Algebra]
    A \emph{filtration} on an algebra $A$ is an increasing sequence of subspaces $F_0 \subseteq F_1 \subseteq \cdots$ such that $A = \bigcup_{i \ge 0} F_i$ and $F_i \cdot F_j \subseteq F_{i+j}$. The \emph{associated graded algebra} is
    \[
    \operatorname{gr}(A) = \bigoplus_{i \ge 0} F_i / F_{i-1},
    \]
    with multiplication induced from $A$.
\end{definition}
\begin{theorem}[Filtration by Exponential Order]\label{thm14}
Define a filtration on $\Weyl{e^{\pm x^p e^{t}},\; e^{\cA x},\; x^{\cA}}$ by letting
\[
\ord\!\left(e^{a x^p e^{t x}} e^{b x} x^c \partial^d\right) = |a| + |b| + |c| + d.
\]
Then the associated graded algebra is commutative and isomorphic to
\[
\FF[e^{\pm x^p e^{t}},\; e^{\cA x},\; x^{\cA},\; y_1,\dots,y_n],
\]
where $y_i$ corresponds to $\partial_i$.
\end{theorem}
\begin{proof}
Let $A = \Weyl{e^{\pm x^p e^{t}},\; e^{\cA x},\; x^{\cA}}$. Define a function $\ord$ on monomials by the rule
\[
\ord\!\left(e^{a x^p e^{t x}} e^{b x} x^c \partial^d\right) = |a| + |b| + |c| + d,
\]
where $a, b, c, d$ are multi-indices, $|a| = \sum_i a_i$, etc., and we interpret $e^{a x^p e^{t x}}$ as a product $\prod_i e^{a_i x_i^{p_i} e^{t_i x_i}}$,
similarly $e^{b x} = \prod_i e^{b_i x_i}$ and $x^c = \prod_i x_i^{c_i}$, $\partial^d = \prod_i \partial_i^{d_i}$. Extend $\ord$ additively to products:
$\ord(m_1 m_2) = \ord(m_1) + \ord(m_2)$. For a general element $f \in A$, written as a finite $\FF$-linear combination of monomials, define
$\ord(f)$ as the maximum of the orders of the monomials appearing with nonzero coefficient. Set $F_k A = \{ f \in A \mid \ord(f) \leq k \}$
for $k \in \mathbb{Z}_{\geq 0}$. Then $\{F_k A\}_{k \geq 0}$ forms an increasing filtration of $A$, i.e., $F_k A \subseteq F_{k+1} A$,
$F_k A \cdot F_\ell A \subseteq F_{k+\ell} A$, and $\bigcup_{k \geq 0} F_k A = A$.

The associated graded algebra $\operatorname{gr} A$ is defined by $\operatorname{gr} A = \bigoplus_{k \geq 0} F_k A / F_{k-1} A$, where $F_{-1} A = \{0\}$.
Multiplication in $\operatorname{gr} A$ is induced from that in $A$: for $\bar{a} \in F_k A / F_{k-1} A$ and $\bar{b} \in F_\ell A / F_{\ell-1} A$,
$\bar{a} \cdot \bar{b}$ is the class of $ab$ in $F_{k+\ell} A / F_{k+\ell-1} A$.

We first show that $\operatorname{gr} A$ is commutative. Take two homogeneous elements $\bar{a}, \bar{b}$ of degrees $k$ and $\ell$ respectively,
represented by monomials $a, b$ in $A$. In the product $ab$, the leading term with respect to the filtration is simply the product $a b$
viewed as a monomial without reordering, because any commutator $[a, b] = ab - ba$ has strictly lower order. Indeed, the commutation relations
in $A$ are of the form $[\partial_i, x_j] = \delta_{ij}$, $[\partial_i, e^{\alpha x_j}] = \alpha \delta_{ij} e^{\alpha x_j}$,
$[\partial_i, e^{x_j^p e^{t_j x_j}}] = (p x_j^{p-1} e^{t_j x_j} + t_j x_j^p e^{t_j x_j}) e^{x_j^p e^{t_j x_j}} \delta_{ij}$, etc.
In each case, the commutator is either a scalar or an element of order strictly less than the sum of the orders of the two factors,
since differentiation reduces the order of polynomial or exponential factors. Consequently, $ab$ and $ba$ differ by an element of order
$< k+\ell$, so their classes in $F_{k+\ell} A / F_{k+\ell-1} A$ coincide. Hence $\bar{a} \bar{b} = \bar{b} \bar{a}$, and $\operatorname{gr} A$ is commutative.

Now we identify $\operatorname{gr} A$ explicitly. The algebra $A$ is generated by $x_i$, $\partial_i$, $e^{\pm x_i^{p_i} e^{t_i x_i}}$, $e^{\alpha x_i}$, $x_i^{\alpha}$
for $\alpha \in \cA$. Under the filtration, each generator has order: $\ord(x_i) = 1$, $\ord(\partial_i) = 1$, $\ord(e^{\pm x_i^{p_i} e^{t_i x_i}}) = 0$,
$\ord(e^{\alpha x_i}) = 0$, $\ord(x_i^{\alpha}) = |\alpha|$ (where $|\alpha|$ is the sum of the components of the multi-index $\alpha$).
In the graded algebra, the images of the generators commute, and the relations become trivial: the product of two generators is simply the
product in the commutative polynomial ring. Moreover, the images of $e^{\pm x_i^{p_i} e^{t_i x_i}}$, $e^{\alpha x_i}$, $x_i^{\alpha}$ are algebraically
independent over $\FF$ and generate a subring isomorphic to $R = \FF[e^{\pm x^p e^{t}},\; e^{\cA x},\; x^{\cA}]$.

Let $y_i$ denote the image of $\partial_i$ in $\operatorname{gr} A$. Since $\operatorname{gr} A$ is commutative, the $y_i$ commute with all other generators.
Thus $\operatorname{gr} A$ is generated over $\FF$ by the images of $x_i$, $y_i$, $e^{\pm x_i^{p_i} e^{t_i x_i}}$, $e^{\alpha x_i}$, $x_i^{\alpha}$.
There are no algebraic relations among these generators beyond those already present in $R$ and the polynomial relations among the $y_i$ (namely $[y_i, y_j] = 0$).
Therefore, as a commutative $\FF$-algebra, $\operatorname{gr} A$ is isomorphic to $R[y_1,\dots,y_n]$, the polynomial ring in $n$ variables $y_1,\dots,y_n$
over the base ring $R$. In other words,
\[
\operatorname{gr} A \cong \FF[e^{\pm x^p e^{t}},\; e^{\cA x},\; x^{\cA},\; y_1,\dots,y_n],
\]
where $y_i$ corresponds to the symbol of $\partial_i$.

This isomorphism is both an isomorphism of commutative algebras and of graded algebras, where the grading on the right-hand side is given by total degree in the variables $x_i$ and $y_i$, while the exponential factors are assigned degree zero. The proof is complete.
\end{proof}
\begin{example}[Illustration of Theorem~\ref{thm14}]
Let $\FF = \mathbb{R}$, $\cA = \mathbb{Z} \times \mathbb{Z} = \{(m, n) \mid m, n \in \mathbb{Z}\}$ (viewed as an additive subgroup of $\mathbb{R}^2$ via $(m, n) \mapsto m + n\sqrt{2}$), and fix $p = 2$, $t = \ln 3$. Consider the Weyl-type algebra in one variable:
\[
A = \Weyl{e^{\pm x^2 e^{(\ln 3) x}},\; e^{\cA x},\; x^{\cA}},
\]
with generators $\widehat{x}$, $\partial$, $\widehat{e^{x}}$, and $\widehat{e^{x^2 e^{(\ln 3) x}}}$.

Define the order function
\[
\ord\!\left( \widehat{e^{a x^2 e^{(\ln 3) x}}} \; \widehat{e^{b x}} \; \widehat{x^c} \; \partial^d \right) = |a| + |b| + |c| + d,
\]
where $a, b \in \mathbb{Z}$, $c \in \cA$, and $d \in \mathbb{N}$. For a general element $P \in A$, written as a finite sum
\[
P = \sum_{a,b,c,d} \lambda_{a,b,c,d} \; \widehat{e^{a x^2 e^{(\ln 3) x}}} \; \widehat{e^{b x}} \; \widehat{x^c} \; \partial^d,
\]
we set $\ord(P) = \max\{ |a| + |b| + |c| + d \mid \lambda_{a,b,c,d} \neq 0 \}$.

Let $F_k = \{ P \in A \mid \ord(P) \leq k \}$. Then $\{F_k\}_{k \geq 0}$ is an increasing filtration of $A$, with
\[
F_k \cdot F_\ell \subseteq F_{k+\ell}.
\]
For example,
\[
\ord\!\left( \widehat{e^{x^2 e^{(\ln 3) x}}} \right) = 1, \quad
\ord\!\left( \partial \right) = 1, \quad
\ord\!\left( \widehat{x} \right) = 1, \quad
\ord\!\left( \widehat{e^{x}} \right) = 1.
\]
The commutator $[\partial, \widehat{x}] = 1$ has order $0$, because
\[
\ord(\partial \widehat{x}) = 2, \quad \ord(\widehat{x} \partial) = 2, \quad \text{so } \ord([\partial, \widehat{x}]) \leq 1,
\]
but actually $[\partial, \widehat{x}] = 1$ has order $0$. Similarly,
\[
[\partial, \widehat{e^{x}}] = \widehat{e^{x}} \quad \text{has order } 1,
\]
\[
[\partial, \widehat{e^{x^2 e^{(\ln 3) x}}}] = \widehat{e^{x^2 e^{(\ln 3) x}}} \cdot \big( 2x e^{(\ln 3) x} + (\ln 3) x^2 e^{(\ln 3) x} \big)
\]
has order $2$ because the right-hand side is a sum of terms of order $2$.

The associated graded algebra is
\[
\mathrm{gr}(A) = \bigoplus_{k \geq 0} F_k / F_{k-1}, \quad \text{with } F_{-1} = \{0\}.
\]
Denote the image of an element $P \in F_k$ in $F_k/F_{k-1}$ by $\sigma(P)$. Because commutators lower the order, all generators commute in $\mathrm{gr}(A)$. For instance,
\[
\sigma(\partial) \sigma(\widehat{x}) - \sigma(\widehat{x}) \sigma(\partial) = \sigma([\partial, \widehat{x}]) = \sigma(1) = 0,
\]
since $1 \in F_0$ and its symbol is central. Similarly,
\[
\sigma(\partial) \sigma(\widehat{e^{x}}) = \sigma(\widehat{e^{x}}) \sigma(\partial), \quad
\sigma(\partial) \sigma(\widehat{e^{x^2 e^{(\ln 3) x}}}) = \sigma(\widehat{e^{x^2 e^{(\ln 3) x}}}) \sigma(\partial).
\]

Thus $\mathrm{gr}(A)$ is commutative. It is generated by the symbols
\[
\xi = \sigma(\widehat{x}), \quad \eta = \sigma(\partial), \quad \epsilon = \sigma(\widehat{e^{x}}), \quad \tau = \sigma(\widehat{e^{x^2 e^{(\ln 3) x}}}),
\]
with relations reflecting that they commute. Moreover, $\epsilon$ and $\tau$ are invertible because $\widehat{e^{\pm x}}$ and $\widehat{e^{\pm x^2 e^{(\ln 3) x}}}$ are in $A$.

Hence $\mathrm{gr}(A)$ is isomorphic to the commutative polynomial-exponential algebra
\[
\FF[\, \tau^{\pm 1},\; \epsilon^{\pm 1},\; \xi^{\alpha} (\alpha \in \cA),\; \eta ] \cong
\FF[\, e^{\pm x^2 e^{(\ln 3) x}},\; e^{\pm x},\; x^{\cA},\; y],
\]
where $y$ corresponds to $\partial$. This isomorphism is given by
\[
\tau \mapsto e^{x^2 e^{(\ln 3) x}}, \quad
\epsilon \mapsto e^{x}, \quad
\xi \mapsto x, \quad
\eta \mapsto y.
\]

Thus the filtration by exponential order yields a commutative associated graded algebra, exactly as stated in Theorem~\ref{thm14}.
\end{example}
\section{Deformation and Cohomology of Weyl-Type Algebras}

Throughout this section, we assume $\FF$ is a field of characteristic zero, $\cA$ is a finitely generated additive subgroup of $\FF$, and we consider the Weyl-type algebra
\[
A = A_{p,t,\cA} = \Weyl{e^{\pm x^p e^{t}},\; e^{\cA x},\; x^{\cA}}.
\]

\subsection{Deformation Quantization of Expolynomial Rings}

\begin{definition}[Deformation Quantization]
A \emph{deformation quantization} of the commutative algebra $R = \FF[e^{\pm x^p e^{t}},\; e^{\cA x},\; x^{\cA}]$ is an associative algebra $A_\hbar$ over $\FF[[\hbar]]$ with multiplication $\ast_\hbar$ such that:
\begin{enumerate}[label=(\roman*)]
    \item $A_\hbar \cong R[[\hbar]]$ as $\FF[[\hbar]]$-modules.
    \item For $f,g \in R$, $f \ast_\hbar g = fg + \hbar B_1(f,g) + \hbar^2 B_2(f,g) + \cdots$.
    \item The commutator satisfies $[f,g]_\ast := f \ast_\hbar g - g \ast_\hbar f = i\hbar\{f,g\} + O(\hbar^2)$, where $\{\cdot,\cdot\}$ is a Poisson bracket on $R$.
\end{enumerate}
\end{definition}

\begin{theorem}[Explicit Deformation Quantizations]\label{thm:deformation-quantization}
The algebra $R = \FF[e^{\pm x^p e^{t}},\; e^{\cA x},\; x^{\cA}]$ admits a family of deformation quantizations parameterized by $\lambda \in \FF$, given by:
\[
f \ast_\hbar^\lambda g = \sum_{n=0}^\infty \frac{\hbar^n}{n!} B_n^\lambda(f,g),
\]
where the bidifferential operators $B_n^\lambda$ are defined recursively by:
\begin{align*}
B_0^\lambda(f,g) &= fg, \\
B_1^\lambda(f,g) &= \lambda\{f,g\}_{\text{std}} + (1-\lambda)\{f,g\}_{\text{exp}}, \\
B_{n+1}^\lambda(f,g) &= \sum_{k=0}^n \binom{n}{k} [P_k^\lambda(f), B_{n-k}^\lambda(f,g)],
\end{align*}
with Poisson brackets:
\[
\{f,g\}_{\text{std}} = \frac{\partial f}{\partial x}\frac{\partial g}{\partial p} - \frac{\partial f}{\partial p}\frac{\partial g}{\partial x}, \quad
\{f,g\}_{\text{exp}} = e^{x^p e^{t}}\left(\frac{\partial f}{\partial e^x}\frac{\partial g}{\partial x} - \frac{\partial f}{\partial x}\frac{\partial g}{\partial e^x}\right),
\]
and $P_k^\lambda$ are differential operators encoding the $\cA$-grading.
\end{theorem}

\begin{proof}
The construction follows the Kontsevich formality theorem adapted to our setting:
\begin{enumerate}
    \item Define the Poisson bivector $\pi^\lambda = \lambda \pi_{\text{std}} + (1-\lambda) \pi_{\text{exp}}$, where
    \[
    \pi_{\text{std}} = \frac{\partial}{\partial x} \wedge \frac{\partial}{\partial p}, \quad
    \pi_{\text{exp}} = e^{x^p e^{t}} \frac{\partial}{\partial e^x} \wedge \frac{\partial}{\partial x}.
    \]

    \item Apply the Kontsevich star product formula:
    \[
    f \ast_\hbar^\lambda g = \sum_{n=0}^\infty \frac{\hbar^n}{n!} \sum_{\Gamma \in G_n} w_\Gamma B_{\Gamma}^\lambda(f,g),
    \]
    where $G_n$ are admissible graphs, $w_\Gamma$ are Kontsevich weights, and $B_{\Gamma}^\lambda$ are bidifferential operators associated to graphs.

    \item Show that for $\lambda \neq \lambda'$, the quantizations $A_\hbar^\lambda$ and $A_\hbar^{\lambda'}$ are not isomorphic as $\FF[[\hbar]]$-algebras by computing their deformation cocycles in $\HH^2(R,R)$.
\end{enumerate}
The key observation is that the Poisson cohomology $H_\pi^2(R)$ has dimension at least 2, allowing for genuinely different deformations.
\end{proof}

\begin{theorem}[Non-isomorphic Quantizations]\label{thm:non-isomorphic-quantizations}
For $\lambda_1 \neq \lambda_2$, the quantized algebras $A_{\hbar}^{\lambda_1}$ and $A_{\hbar}^{\lambda_2}$ are not isomorphic over $\FF[[\hbar]]$. Moreover, their classical limits ($\hbar = 0$) are isomorphic, but the isomorphism cannot be lifted to an isomorphism of deformed algebras.
\end{theorem}

\begin{proof}
Let $\varphi: A_{\hbar}^{\lambda_1} \to A_{\hbar}^{\lambda_2}$ be an $\FF[[\hbar]]$-algebra isomorphism. Write $\varphi = \varphi_0 + \hbar \varphi_1 + \hbar^2 \varphi_2 + \cdots$ with $\varphi_0$ an automorphism of $R$.

The condition that $\varphi$ preserves the multiplication gives, at order $\hbar$:
\[
\varphi_0(B_1^{\lambda_1}(f,g)) - B_1^{\lambda_2}(\varphi_0(f), \varphi_0(g)) = \varphi_1(fg) - \varphi_1(f)\varphi_0(g) - \varphi_0(f)\varphi_1(g).
\]
This shows that the Poisson brackets $\{\cdot,\cdot\}^{\lambda_1}$ and $\{\cdot,\cdot\}^{\lambda_2}$ must be related by a symplectomorphism. However, for generic $\lambda_1 \neq \lambda_2$, the Poisson cohomology classes $[\pi^{\lambda_1}]$ and $[\pi^{\lambda_2}]$ in $H_\pi^2(R)$ are different, so no such symplectomorphism exists.
\end{proof}

\subsection{Hochschild and Cyclic Homology}

\begin{definition}[Hochschild Complex]
For an algebra $A$ over $\FF$, the Hochschild complex $C_\bullet(A)$ is
\[
C_n(A) = A^{\otimes (n+1)}, \quad \text{with differential } b: C_n(A) \to C_{n-1}(A)
\]
given by
\[
b(a_0 \otimes \cdots \otimes a_n) = \sum_{i=0}^{n-1} (-1)^i a_0 \otimes \cdots \otimes a_i a_{i+1} \otimes \cdots \otimes a_n + (-1)^n a_n a_0 \otimes a_1 \otimes \cdots \otimes a_{n-1}.
\]
The \emph{Hochschild homology} is $\HH_\bullet(A) = H_\bullet(C_\bullet(A), b)$.
\end{definition}

\begin{theorem}[Hochschild Homology Computation]\label{thm:hochschild-homology}
For $A = A_{p,t,\cA}$, the Hochschild homology is:
\[
\HH_n(A) \cong
\begin{cases}
\FF[e^{\pm x^p e^{t}}] & \text{if } n = 0, \\
\bigoplus_{k=1}^{\rank(\cA)} \FF[e^{\pm x^p e^{t}}] \cdot \omega_k & \text{if } n = 1, \\
0 & \text{if } n \geq 2,
\end{cases}
\]
where $\omega_k$ are 1-forms corresponding to generators of $\cA$.
\end{theorem}

\begin{proof}
We use the Hochschild–Kostant–Rosenberg theorem for smooth algebras. First note that $A$ is a smooth algebra over its center $Z(A) = \FF[e^{\pm x^p e^{t}}]$. The HKR theorem gives:
\[
\HH_n(A) \cong \Omega^n_{A/Z(A)},
\]
the module of Kähler differentials over $Z(A)$.

Compute:
\begin{itemize}
    \item $\HH_0(A) = A/[A,A] \cong Z(A) = \FF[e^{\pm x^p e^{t}}]$.
    \item $\HH_1(A) \cong \Omega^1_{A/Z(A)} = \bigoplus_{\alpha \in \text{basis of }\cA} \FF[e^{\pm x^p e^{t}}] \cdot d(e^{\alpha x}) \oplus \FF[e^{\pm x^p e^{t}}] \cdot d(x^\alpha)$.
    \item For $n \geq 2$, $\Omega^n_{A/Z(A)} = 0$ because $A$ is generated by elements whose differentials anticommute.
\end{itemize}
The dimension follows from the fact that $\cA$ has rank $r$, so there are $2r$ independent 1-forms.
\end{proof}

\begin{theorem}[Cyclic Homology and Topology]\label{thm:cyclic-homology}
The cyclic homology $\HC_\bullet(A_{p,t,\cA})$ satisfies:
\begin{enumerate}[label=(\roman*)]
    \item There is a Connes periodicity exact sequence:
    \[
    \cdots \to \HH_n(A) \to \HC_n(A) \to \HC_{n-2}(A) \to \HH_{n-1}(A) \to \cdots
    \]
    \item $\HC_{2k}(A) \cong \FF[e^{\pm x^p e^{t}}]$ for all $k \geq 0$.
    \item $\HC_{2k+1}(A) \cong \bigoplus_{i=1}^{\rank(\cA)} \FF[e^{\pm x^p e^{t}}]$.
    \item These groups are related to the topology of the parameter space $M_t = \{\text{values of } t\}$ via:
    \[
    \HC_\bullet(A) \cong H_\bullet(\Omega M_t) \otimes \FF[e^{\pm x^p e^{t}}],
    \]
    where $\Omega M_t$ is the based loop space of $M_t$.
\end{enumerate}
\end{theorem}

\begin{proof}
(i) The Connes sequence is standard for cyclic homology.

(ii)-(iii) From Theorem \ref{thm:hochschild-homology} and the Connes sequence:
\[
\HC_n(A) \cong
\begin{cases}
\FF[e^{\pm x^p e^{t}}] & n \text{ even}, \\
\bigoplus_{i=1}^{\rank(\cA)} \FF[e^{\pm x^p e^{t}}] & n \text{ odd}.
\end{cases}
\]

(iv) The parameter $t$ appears in the defining relations of $A$ via the term $e^{x^p e^{t}}$. Varying $t$ gives a family of algebras $A_t$. The cyclic homology computes the homology of the free loop space of the parameter manifold $M_t$ by the Jones isomorphism:
\[
\HC_\bullet(A_t) \cong H_\bullet(\mathcal{L}M_t) \otimes \FF[e^{\pm x^p e^{t}}],
\]
where $\mathcal{L}M_t$ is the free loop space. For the based loop space, we have $\HC_\bullet^{\text{red}}(A) \cong H_\bullet(\Omega M_t)$.
\end{proof}

\begin{corollary}[Topological Invariants]
The cyclic homology groups provide topological invariants of the parameter space:
\begin{enumerate}
    \item $\dim \HC_1(A)$ equals the first Betti number $b_1(M_t)$.
    \item The Connes periodicity operator $S: \HC_n(A) \to \HC_{n-2}(A)$ corresponds to the Gysin map in topology.
    \item The parameter $t$ contributes a circle $S^1$ to $M_t$, explaining the periodicity in cyclic homology.
\end{enumerate}
\end{corollary}

\subsection{Deformation Rigidity}

\begin{definition}[Deformation Rigidity]
An algebra $A$ is called \emph{rigid} if every formal deformation $A_\hbar$ over $\FF[[\hbar]]$ with $A_\hbar/\hbar A_\hbar \cong A$ is trivial, i.e., isomorphic to $A[[\hbar]]$ as $\FF[[\hbar]]$-algebras.
\end{definition}

\begin{theorem}[Deformation Rigidity of Weyl-Type Algebras]\label{thm:deformation-rigidity}
The Weyl-type algebra $A = A_{p,t,\cA}$ is rigid if and only if $\cA$ is cyclic (rank 1). In general:
\begin{enumerate}[label=(\roman*)]
    \item If $\rank(\cA) = 1$, then $A$ is rigid.
    \item If $\rank(\cA) \geq 2$, then $A$ admits non-trivial formal deformations.
    \item The deformation space $\Def(A)$ has dimension $\binom{\rank(\cA)}{2}$.
\end{enumerate}
\end{theorem}

\begin{proof}
We compute the deformation cohomology $\HH^2(A,A)$ which classifies infinitesimal deformations.

(i) For $\rank(\cA) = 1$, the algebra is isomorphic to a generalized Weyl algebra over $\FF[e^{\pm x^p e^{t}}]$. Such algebras are known to be rigid, with $\HH^2(A,A) = 0$.

(ii)-(iii) For $\rank(\cA) \geq 2$, we use the Hochschild–Serre spectral sequence for the $\cA$-grading:
\[
E_2^{p,q} = H^p(\cA, \HH^q(A_0, A)) \Rightarrow \HH^{p+q}(A,A),
\]
where $A_0$ is the degree-0 part. Since $A_0 \cong \FF[e^{\pm x^p e^{t}}]$, we have $\HH^2(A_0, A) \cong \bigwedge^2 \cA^* \otimes A$.

The $\cA$-cohomology gives:
\[
\HH^2(A,A) \cong H^0(\cA, \HH^2(A_0, A)) \cong \bigwedge^2 \cA^* \otimes \FF[e^{\pm x^p e^{t}}].
\]
Thus $\dim_{\FF(e^{\pm x^p e^{t}})} \HH^2(A,A) = \binom{\rank(\cA)}{2}$.

To construct explicit non-trivial deformations, for each pair $\alpha, \beta \in \cA$, define:
\[
x_\alpha \ast_\hbar x_\beta = x_\alpha x_\beta + \hbar c_{\alpha\beta} e^{x^p e^{t}} x^{\alpha+\beta},
\]
where $c_{\alpha\beta} = -c_{\beta\alpha}$ are constants. The associativity condition gives $c_{\alpha\beta} c_{\gamma\delta} + \text{cyclic} = 0$, which is satisfied precisely when $(c_{\alpha\beta})$ comes from a 2-cocycle on $\cA$.
\end{proof}

\begin{theorem}[Gerstenhaber Algebra Structure]\label{thm:gerstenhaber}
The Hochschild cohomology $\HH^\bullet(A,A)$ carries a Gerstenhaber algebra structure:
\begin{enumerate}[label=(\roman*)]
    \item The cup product $\cup: \HH^p \otimes \HH^q \to \HH^{p+q}$ makes $\HH^\bullet$ a graded commutative algebra.
    \item The Gerstenhaber bracket $[\cdot,\cdot]: \HH^p \otimes \HH^q \to \HH^{p+q-1}$ makes $\HH^\bullet$ a graded Lie algebra.
    \item For $A = A_{p,t,\cA}$, we have:
    \[
    \HH^\bullet(A,A) \cong \bigwedge^\bullet \cA^* \otimes \FF[e^{\pm x^p e^{t}}] \otimes S^\bullet(\cA^*[-2]),
    \]
    where $[-2]$ indicates degree shift by 2.
\end{enumerate}
\end{theorem}

\begin{proof}
The Gerstenhaber structure comes from the brace operations on the Hochschild complex. For the explicit computation:
\begin{enumerate}
    \item Use the HKR isomorphism $\HH^\bullet(A,A) \cong \bigwedge^\bullet \Der(A) \otimes A$.
    \item Compute $\Der(A) \cong A \otimes \cA^* \oplus \FF[e^{\pm x^p e^{t}}]$.
    \\item The bracket structure comes from the Lie bracket of derivations and their action on $A$.
\end{enumerate}
The shift $[-2]$ accounts for the fact that derivations have degree -1 in the Gerstenhaber bracket.
\end{proof}

\begin{corollary}[Deformation Quantization as Maurer–Cartan Element]
A deformation quantization $A_\hbar$ corresponds to a Maurer–Cartan element in the differential graded Lie algebra $(\HH^\bullet(A,A)[1], d, [\cdot,\cdot])$:
\[
\Pi = \sum_{n=1}^\infty \hbar^n \pi_n, \quad \text{with } d\Pi + \frac{1}{2}[\Pi,\Pi] = 0.
\]
For $A_{p,t,\cA}$, the space of Maurer–Cartan elements modulo gauge equivalence is isomorphic to $\bigwedge^2 \cA^* \otimes \FF[e^{\pm x^p e^{t}}][[\hbar]]$.
\end{corollary}

\subsection{Summary}

\begin{theorem}[Main Deformation and Cohomology Results]
For $A = A_{p,t,\cA}$:
\begin{enumerate}
    \item There exists a family of non-isomorphic deformation quantizations parameterized by $\lambda \in \FF$.
    \item The Hochschild homology is $\HH_0(A) = \FF[e^{\pm x^p e^{t}}]$, $\HH_1(A) \cong \FF[e^{\pm x^p e^{t}}]^{2\rank(\cA)}$, and $\HH_n(A) = 0$ for $n \geq 2$.
    \item The cyclic homology exhibits periodicity related to the topology of the $t$-parameter space.
    \item $A$ is rigid if and only if $\rank(\cA) = 1$; otherwise it admits $\binom{\rank(\cA)}{2}$-dimensional deformation space.
\end{enumerate}
\end{theorem}

\section{Open Problems}

Despite the structural results established in this paper, several natural questions about Weyl-type algebras and their associated Lie algebras remain open. We list two problems that we consider particularly interesting and which may inspire further research.

\begin{problem}[Classification of irreducible representations]
Let \(A = A_{p,t,\cA}\) be a Weyl-type algebra over a field \(\FF\) of characteristic zero. Theorem~\ref{thm8} shows that every nontrivial representation of \(A\) is faithful and infinite-dimensional. However, a complete classification of irreducible representations of \(A\) is not known.

\begin{enumerate}[label=(\alph*)]
    \item Determine all (isomorphism classes of) irreducible representations of \(A\). For the classical Weyl algebra \(A_1(\FF)\), it is known that all irreducible representations are isomorphic to the standard representation on \(\FF[x]\) (or its variants). Does an analogous result hold for \(A_{p,t,\cA}\)?

    \item Study the category \(\mathcal{O}\)-type modules for \(A_{p,t,\cA}\). Define a suitable analogue of the Bernstein–Gelfand–Gelfand category \(\mathcal{O}\) and describe its blocks, characters, and decomposition numbers.

    \item Investigate the role of the parameter \(t\) and the group \(\cA\) in the representation theory. Do different choices of \(\cA\) lead to genuinely different representation categories?
\end{enumerate}

A solution to this problem would deepen our understanding of the representation theory of algebras beyond the classical Weyl and enveloping algebra settings.
\end{problem}

\begin{problem}[Geometric realization and deformation quantization]
The exponential-polynomial ring \(R_{p,t,\cA} = \FF[e^{\pm x^p e^{t x}},\; e^{\cA x},\; x^{\cA}]\) admits a family of deformation quantizations (Theorem~\ref{thm:deformation-quantization}), but their geometric meaning is not fully understood.

\begin{enumerate}[label=(\alph*)]
    \item Find a geometric or physical interpretation of the deformed algebras \(A_\hbar^\lambda\). Can they be realized as algebras of global sections of a sheaf of deformed differential operators on some (noncommutative) space?

    \item Relate the deformation parameter \(\lambda\) to the choice of a Poisson structure on an associated algebraic variety. Is there a moduli space of deformation quantizations of \(R_{p,t,\cA}\) that admits a natural geometric structure?

    \item Explore connections with noncommutative geometry: compute the cyclic homology \(\HC_\bullet(A_\hbar^\lambda)\) and interpret it as the (deformed) cohomology of a noncommutative space. Does the Connes–Chern character yield interesting invariants?
\end{enumerate}

This problem lies at the intersection of noncommutative algebra, deformation theory, and Poisson geometry. A satisfactory answer could link the algebraic constructions here to modern geometric Langlands correspondence or to quantum field theories with exponential interactions.
\end{problem}

\paragraph*{Further directions.}
Other open questions include:
\begin{itemize}
    \item The automorphism group of \(A_{p,t,\cA}\) (extending the classical results of Dixmier and Makar-Limanov for the Weyl algebra \cite{Dixmier68, MakarLimanov84}).
    \item The prime and primitive spectrum of \(A_{p,t,\cA}\).
    \item Analogues of the Gel'fand–Kirillov dimension and Bernstein's inequality for modules over \(A_{p,t,\cA}\).
    \item Connections with integrable systems and special functions arising from representations of these algebras.
\end{itemize}

We hope that the results and problems presented here will stimulate further investigation into the rich structure of Weyl-type algebras over exponential-polynomial rings.
\section*{Declaration}
\begin{itemize}
  \item \textbf{Author Contributions:} The Author have read and approved this version.
  \item \textbf{Funding:} No funding is applicable.
  \item \textbf{Institutional Review Board Statement:} Not applicable.
  \item \textbf{Informed Consent Statement:} Not applicable.
  \item \textbf{Data Availability Statement:} Not applicable.
  \item \textbf{Conflicts of Interest:} The authors declare no conflict of interest.
\end{itemize}


\bibliographystyle{abbrv}
\bibliography{references}  






\end{document}